\documentclass[11pt,leqno,twoside,a4paper]{article}
\usepackage[english]{babel}
\usepackage[utf8]{inputenc}
\usepackage{authblk}
\usepackage{amsfonts,amssymb,amsmath,amsthm}
\usepackage{enumitem}
\usepackage[a4paper,top=25mm,bottom=30mm,left=30mm,right=30mm]{geometry}
\usepackage{graphicx} 

\usepackage{verbatim}
\usepackage{calc}
\usepackage{color}
\usepackage{xcolor}
\usepackage{lineno}
\usepackage[mathscr]{eucal}
\usepackage{needspace}
\usepackage{lmodern}     
\usepackage[T1]{fontenc}
\usepackage{microtype}
\usepackage{ifthen}
\usepackage{esint,braket,mathtools,comment}
\usepackage[hidelinks,hyperfootnotes=false,hypertexnames=false]{hyperref}
\usepackage{doi}
\usepackage{xurl}
\renewcommand{\doi}[1]{\doitext\href{https://doi.org/#1}{\nolinkurl{#1}}}
\usepackage[capitalise]{cleveref}

\providecommand{\appendixname}{Appendix}

\usepackage{fancyhdr}
\newlength{\abstractleft}
\newlength{\abstractright}
\makeatletter
\renewenvironment{abstract}{%
  \par\small
  \list{}{%
    \setlength{\leftmargin}{\abstractleft}%
    \setlength{\rightmargin}{\abstractright}%
    \setlength{\listparindent}{0pt}%
    \setlength{\itemindent}{0pt}%
    \setlength{\parsep}{0pt}%
  }%
\item\relax
  \noindent\textsc{Abstract.}\ %
  \ignorespaces
}{%
  \endlist
  \par\normalsize
}
\makeatother



\newcommand{\ShortAuthors}{\scriptsize H. CASTRO AND A. ZUNIGA} 
\newcommand{\ShortTitle}{\scriptsize
Sharp magnetic Hardy--trace inequalities} 

\fancypagestyle{plain}{%
  \fancyhf{}
  \fancyfoot[C]{\thepage}%
}

\usepackage{etoolbox}
\newif\ifinappendix
\pretocmd{\appendix}{%
  \inappendixtrue
  \addtocontents{toc}{\protect\inappendixtrue}%
}{}{}

\makeatletter

\renewcommand{\@seccntformat}[1]{%
  \ifinappendix
  \ifdefstring{#1}{section}{\appendixname~\thesection.\quad}{%
    \ifdefstring{#1}{subsection}{\appendixname~\thesubsection.\quad}{%
      \csname the#1\endcsname\quad
    }%
  }%
  \else
  \csname the#1\endcsname\quad
  \fi
}
\makeatletter

\newcommand{\TOCEntryFont}{}

\renewcommand{\l@section}[2]{%
  \begingroup
  \parskip=0pt
  \let\oldnumberline\numberline
  \ifinappendix
  \renewcommand{\numberline}[1]{\oldnumberline{\appendixname~##1.}}%
  \@dottedtocline{1}{1.5em}{6.8em}{{\TOCEntryFont #1}}{{\TOCEntryFont #2}}%
  \else
  \renewcommand{\numberline}[1]{\oldnumberline{##1.}}%
  \@dottedtocline{1}{1.5em}{3.0em}{{\TOCEntryFont #1}}{{\TOCEntryFont #2}}%
  \fi
  \endgroup
}

\renewcommand{\l@subsection}[2]{%
  \begingroup
  \parskip=0pt
  \let\oldnumberline\numberline
  \ifinappendix
  \renewcommand{\numberline}[1]{\oldnumberline{\appendixname~##1.}}%
  \@dottedtocline{2}{3.8em}{7.4em}{{\TOCEntryFont #1}}{{\TOCEntryFont #2}}%
  \else
  \renewcommand{\numberline}[1]{\oldnumberline{##1.}}%
  \@dottedtocline{2}{3.8em}{3.8em}{{\TOCEntryFont #1}}{{\TOCEntryFont #2}}%
  \fi
  \endgroup
}
\makeatother

\usepackage{titlesec}
\titleformat{\section}
{\normalfont\scshape\large\filcenter}
{\ifinappendix\appendixname~\thesection.\else\thesection.\fi}{0.8em}{}

\titleformat{\subsection}
{\normalfont\scshape\normalsize\filcenter}
{\ifinappendix\appendixname~\thesubsection.\else\thesubsection.\fi}{0.8em}{}

\titlespacing*{\section}{0pt}{1.5ex plus .2ex}{0.8ex}
\titlespacing*{\subsection}{0pt}{1.2ex plus .2ex}{0.6ex}

\titleformat{\subsubsection}
{\normalfont\scshape\normalsize}
{\thesubsubsection.}{0.8em}{}

\newtheorem{theorem}{Theorem}[section]
\newtheorem{lemma}[theorem]{Lemma}
\newtheorem{proposition}[theorem]{Proposition}

\newtheorem{remark}[theorem]{Remark}

\allowdisplaybreaks

\renewcommand{\a}{{\boldsymbol{a}}}
\newcommand{\A}{\mathbf{A}}
\newcommand{\Aext}{\mathbb{A}}
\newcommand{\R}{\mathbb{R}}
\newcommand{\Rr}{\mathcal{R}_{\beta_1,\beta_2}}
\newcommand{\C}{\mathbb{C}}
\newcommand{\Z}{\mathbb{Z}}

\newcommand{\Ss}{\mathbb{S}}
\newcommand{\Sph}{\Ss^{N-1}}
\newcommand{\ii}{\mathrm{i}}
\newcommand{\dd}{\,\mathrm{d}}
\newcommand{\Kstar}{K^*_{N,\a,\alpha,\beta_1,\beta_2}}
\newcommand{\F}{\,{}_2F_1}
\newcommand{\Vadm}{V_{\mathrm{adm}}}
\newcommand{\lam}{\lambda_{1,\a}}
\newcommand{\ve}{\varepsilon}

\DeclareMathOperator{\diver}{div}
\DeclareMathOperator{\dist}{dist}

\renewcommand{\Re}{\operatorname{Re}}

\newcommand{\pt}[1]{\left(#1\right)}

\newcommand{\abs}[1]{\left|#1\right|}
\newcommand{\norm}[1]{\left\|#1\right\|}

\providecommand{\keywords}[1]
{
  \small
  \textit{Keywords: #1}
}

\title{\large\bfseries\MakeUppercase{
    Sharp weighted Hardy--trace inequalities for homogeneous magnetic potentials
}}
\author{\small\MakeUppercase{Hern\'an Castro\; and\; Andres Zuniga}\textsuperscript{*}}
\date{}

\newcommand{\addressblock}[1]{%
  \par\noindent
  {\scshape #1}\par\medskip
}

\newcommand{\EmailB}[1]{%
  \textit{E-mail address}: %
  {\upshape\ttfamily\bfseries\href{mailto:#1}{#1}}%
}

\makeatletter
\renewcommand{\tableofcontents}{%
  \@starttoc{toc}%
}
\makeatother

\hypersetup{
  pdftitle={Sharp weighted Hardy--trace inequalities for homogeneous magnetic potentials},
  pdfauthor={Hernan Castro and Andres Zuniga},
  pdfsubject={Sharp magnetic Hardy--trace inequalities and their deficit}
}

\begin{document}
\maketitle
\thispagestyle{plain}
\vspace{-2.25em}

\begingroup
\renewcommand{\thefootnote}{\fnsymbol{footnote}}
\setcounter{footnote}{1}
\footnotetext{\raggedright
  \textit{Corresponding author:} Andres Zuniga
  (\href{mailto:andres.zuniga@uoh.cl}{\texttt{andres.zuniga@uoh.cl}}).\par
  \textit{Date:} November 02, 2026.\par
  \par
  \par
  \textit{Affiliations:} Hern\'an Castro, Instituto de Matem\'aticas, Universidad de Talca, Chile;
  Andres Zuniga, Instituto de Ciencias de la Ingenier\'ia, Universidad de O'Higgins, Chile.
}
\endgroup

\setcounter{footnote}{0}

\begin{abstract}
  We prove a sharp magnetic improvement of the family of weighted inequalities on the upper half-space \(\R^{N+1}_+\) that interpolate between the Hardy inequality and the trace Hardy inequality. The magnetic potentials are of Aharonov--Bohm type, \(\A(x)=|x|^{-1}\a(x/|x|)\), where the tangential field \(\a\) may be singular on a closed subset of the sphere. We compute the sharp constant explicitly: it is given by the constant of the corresponding nonmagnetic inequality, evaluated in a virtual dimension that depends on the first eigenvalue of the magnetic Laplacian on the sphere. In particular, the magnetic trace Hardy constant strictly exceeds its Euclidean value whenever this eigenvalue is positive. The proof rests on a ground state representation of the deficit, which also shows that the sharp constant is not attained in the energy space.

  As a consequence, we obtain the sharp Hardy inequality for the fractional powers \((-\Delta_\A)^\sigma\), \(0<\sigma<1\), of the magnetic Laplacian: the constant is the same on smooth functions supported away from the singular set and on the form domain of \((-\Delta_\A)^\sigma\), and it is not attained. Examples include the planar Aharonov--Bohm potential with nonintegral flux and, in dimension three, the Aharonov--Bohm potential and the Dirac monopole.
\end{abstract}

\begin{center}
  {\large\textsc{Contents}}
\end{center}
\vspace{-.75em}
\setcounter{tocdepth}{1}
\tableofcontents
\numberwithin{equation}{section}
\vspace{1em}

\section{Introduction}

The classical Hardy inequality states that, for every \(N\geq1\),
\begin{equation}\label{eq:hardy}
  \pt{\frac{N-2}{2}}^2\int_{\R^N}\frac{\abs{f(x)}^2}{\abs{x}^2}\dd x
  \leq
  \int_{\R^N}\abs{\nabla f(x)}^2\dd x,
  \qquad
  f\in C_c^\infty(\R^N\setminus\{0\}),
\end{equation}
and that the constant \(\pt{(N-2)/2}^2\) is sharp and not attained in the energy space. Observe that for \(N=2\) this constant vanishes; this is no coincidence as the infimum of the Rayleigh quotient is zero:
\[
  \inf_{f\neq0}\frac{\displaystyle\int_{\R^2}\abs{\nabla f(x)}^2\dd x}{\displaystyle\int_{\R^2}\frac{\abs{f(x)}^2}{\abs{x}^2}\dd x}=0.
\]
Nevertheless, in their study of Hardy inequalities for magnetic Dirichlet forms, Laptev and Weidl \cite{LaWe1999} discovered that a nontrivial Hardy inequality in the plane can be recovered in the presence of a suitable magnetic field. To state their result, let \(\A:\R^2\setminus\{0\}\to\R^2\) be the Aharonov--Bohm potential
\begin{equation}\label{eq:AB-potential}
  \A(x)=\bar\a\,\frac{(-x_2,x_1)}{\abs{x}^2},
  \qquad x=(x_1,x_2)\in\R^2\setminus\{0\},
  \qquad \bar\a\in\R,
\end{equation}
where \(\bar\a\) is the normalized flux of the corresponding magnetic field, and let \(\nabla_\A:=\nabla-\ii\A\) be the associated magnetic gradient. Physically, this potential models an infinitely thin solenoid: its enclosed magnetic flux affects quantum interference even though the magnetic field vanishes outside the solenoid \cite{AB1959}. Laptev and Weidl \cite[Theorem 3]{LaWe1999} proved, among other things, that
\begin{equation}\label{eq:hardy-magnetic2}
  \dist(\bar\a,\Z)^2\int_{\R^2}\frac{\abs{f(x)}^2}{\abs{x}^2}\dd x
  \leq
  \int_{\R^2}\abs{\nabla_\A f(x)}^2\dd x,
  \qquad f\in C_c^\infty(\R^2\setminus\{0\}),
\end{equation}
and that this constant is sharp. Balinsky \cite{Bal2003} extended this inequality to Aharonov--Bohm potentials with several singularities.

For a general magnetic potential, the basic information comes from the diamagnetic inequality (see for instance the book of Lieb and Loss \cite{LiebLoss1997}): if \(\A\) is a smooth real vector field in \(\R^N\) and \(\nabla_\A:=\nabla-\ii\A\), then
\begin{equation}\label{eq:diamagnetic}
  \abs{\nabla \abs{f}(x)} \leq \abs{\nabla_\A f(x)}
  \quad\text{for a.e. }x\in\R^N,
  \qquad f\in C^\infty(\R^N),
\end{equation}
where \(|f|\in H^1_{\mathrm{loc}}(\R^N)\), and its gradient is understood weakly. Applying \eqref{eq:hardy} to \(\abs{f}\) and using \eqref{eq:diamagnetic} gives
\begin{equation}\label{eq:hardy-trivial}
  \pt{\frac{N-2}{2}}^2\int_{\R^N}\frac{\abs{f(x)}^2}{\abs{x}^2}\dd x
  \leq
  \int_{\R^N}\abs{\nabla_\A f(x)}^2\dd x,
  \qquad f\in C_c^\infty(\R^N\setminus\{0\}).
\end{equation}
This bound holds for every potential, but its constant does not show the possible dependence on the magnetic field.

To see the effect of the magnetic field, we consider potentials that, like \eqref{eq:AB-potential}, are tangential and homogeneous of degree \(-1\), and we allow them to be singular along rays through the origin. More precisely, let \(\Sigma\subset\Sph\) be a closed set of zero surface measure, and let \(\a:\Sph\setminus\Sigma\to T\Sph\) be a tangential vector field that is smooth on \(\Sph\setminus\Sigma\). The set \(\Sigma\) is introduced so that our framework includes fields that are singular at some points of the sphere, such as the three-dimensional Aharonov--Bohm potential and the Dirac monopole of \cref{sec:example-AB3D,sec:example-monopole}. The Dirac monopole represents a hypothetical isolated magnetic charge, whose existence would require electric charge to be quantized \cite{Dirac1931}. On \(\Sph\setminus\Sigma\), the magnetic gradient is \(\nabla_\a:=\nabla_{\Sph}-\ii\a\), and the magnetic Laplace--Beltrami operator acts as
\begin{equation}\label{eq:LA}
  -\Delta_\a \psi :=
  -\Delta_{\Sph}\psi
  +\ii\diver_{\Sph}(\a)\psi
  +2\ii \a\cdot \nabla_{\Sph}\psi+
  |\a|^2 \psi.
\end{equation}
When \(\Sigma\) is nonempty, this differential expression must be supplemented by a choice of behavior near the excluded set. We make this choice by closing the quadratic form of \((-\Delta_\a)\) on \(C_c^\infty(\Sph\setminus\Sigma)\): from now on, \((-\Delta_\a)\) denotes the Friedrichs extension of its restriction to \(C_c^\infty(\Sph\setminus\Sigma)\), whose form domain is
\begin{equation}\label{eq:HA1}
  H_\a^1(\Sph):=
  \overline{C_c^\infty(\Sph\setminus\Sigma)}^{\norm{\cdot}_{H_\a^1}},
  \qquad
  \norm{\psi}_{H_\a^1}^2:=\int_{\Sph}\bigl(|\nabla_\a\psi|^2+|\psi|^2\bigr)\dd\omega,
\end{equation}
and its first eigenvalue is
\begin{equation}\label{eq:lambda1}
  \lam
  :=
  \inf_{\psi\in H_\a^1(\Sph)\setminus\{0\}}
  \frac{\int_{\Sph} \abs{\nabla_\a\psi}^2\dd\omega}{\int_{\Sph} \abs{\psi}^2\dd\omega}.
\end{equation}
Throughout, \(\Sigma\) is fixed, and we do not display it in the notation for \(\lam\) and for the optimal constants. Let
\[
  \Gamma_\Sigma:=\{0\}\cup\{r\omega\in\R^N:\ r>0,\ \omega\in\Sigma\}
\]
be the cone over \(\Sigma\) together with its vertex, so that \(\Gamma_\Sigma=\{0\}\) when \(\Sigma=\emptyset\). We consider magnetic potentials \(\A:\R^N\setminus\Gamma_\Sigma\to\R^N\) of Aharonov--Bohm type,
\[
  \A(x)=\frac{1}{\abs{x}}\,\a\pt{\frac{x}{\abs{x}}}.
\]
Such a potential is transverse, \(\A(x)\cdot x=0\), and smooth away from \(\Gamma_\Sigma\); observe that the potential \eqref{eq:AB-potential} corresponds to \(N=2\), \(\Sigma=\emptyset\) and \(\a(\omega)=\bar\a(-\omega_2,\omega_1)\) for \(\omega=(\omega_1,\omega_2)\in\Ss^1\). Since \(\A\) does not affect the radial derivative, the one-dimensional Hardy inequality in the radial variable, together with the definition \eqref{eq:lambda1} of \(\lam\) applied on each sphere centered at the origin, gives, for every \(f\in C_c^\infty(\R^N\setminus\Gamma_\Sigma)\),
\begin{equation}\label{eq:hardy-magnetic}
  \left(\pt{\frac{N-2}{2}}^2+\lam\right)
  \int_{\R^N} \frac{\abs{f(x)}^2}{\abs{x}^2} \dd x
  \leq
  \int_{\R^N} \abs{\nabla_\A f(x)}^2 \dd x.
\end{equation}
Potentials of this form, with \(\Sigma=\emptyset\), were studied by Felli, Ferrero and Terracini \cite{FeFeTe2011} in their analysis of the behavior of solutions near an isolated singularity of the electromagnetic potential; they showed that, for smooth angular coefficients, the constant in \eqref{eq:hardy-magnetic} is optimal \cite[Lemma 2.1]{FeFeTe2011}. Inequality \eqref{eq:hardy-magnetic} improves the diamagnetic bound \eqref{eq:hardy-trivial} by the first eigenvalue \(\lam\), and it generalizes the result of Laptev and Weidl: for the potential \eqref{eq:AB-potential} one has \(\lam=\dist(\bar\a,\Z)^2\) (see \cref{sec:example-2D}), so that \eqref{eq:hardy-magnetic} reduces to \eqref{eq:hardy-magnetic2}.

The magnetic Hardy inequality \eqref{eq:hardy-magnetic} fits into a broader family of weighted inequalities on the upper half-space, which interpolate between Hardy inequalities and trace inequalities. Let
\[
  \R^{N+1}_+:=\R^N\times(0,\infty),
  \qquad
  \Aext(x,y):=(\A(x),0),
  \qquad
  \nabla_\Aext:=\nabla_{x,y}-\ii\Aext,
\]
and let \(\alpha\in(-1,1)\). Recall the weighted Hardy inequality, as stated by Tzirakis \cite[(11)]{Tzi2015},
\[
  \frac{(N+\alpha-1)^2}{4}\int_{\R^{N+1}_+}\frac{y^\alpha \abs{u}^2}{\abs{x}^2+y^2}\dd x\dd y
  \leq
  \int_{\R^{N+1}_+} y^\alpha \abs{\nabla u}^2\dd x\dd y.
\]
This inequality can be proved using only the derivative with respect to \(\abs{(x,y)}\), so the angular energy can be estimated separately through \eqref{eq:lambda1}. This gives the magnetic version
\begin{equation}\label{sharp-hardy-magnetic-inequality}
  \begin{aligned}
    \frac{(N+\alpha-1)^2}{4}\int_{\R^{N+1}_+}
    \frac{y^\alpha|u|^2}{|x|^2+y^2}\dd x\dd y
    &+\lam\int_{\R^{N+1}_+}\frac{y^\alpha|u|^2}{|x|^2}\dd x\dd y\\
    &\hspace{4em}\leq\int_{\R^{N+1}_+}y^\alpha|\nabla_\Aext u|^2\dd x\dd y,
  \end{aligned}
\end{equation}
whose first coefficient is sharp for the given second coefficient. This sharpness is proved, among other results on inequalities of this type, in the companion paper \cite{CZB}, which also determines the exact region of admissible pairs of coefficients and obtains remainder terms involving the higher eigenvalues of \((-\Delta_\a)\).

At the other end of this family is the trace Hardy inequality, also known as Kato's inequality, whose sharp constant goes back to Herbst \cite{Her1977}. Its weighted version reads
\[
  H_{N,\alpha}\int_{\R^N}\frac{|u(x,0)|^2}{|x|^{1-\alpha}}\dd x
  \leq\int_{\R^{N+1}_+}y^\alpha|\nabla u|^2\dd x\dd y,
\]
and the sharp constant is
\[
  H_{N,\alpha}=
  2\frac{\Gamma\pt{\frac{1+\alpha}{2}}}{\Gamma\pt{\frac{1-\alpha}{2}}}
  \frac{\Gamma\pt{\frac{N+1-\alpha}{4}}^2}{\Gamma\pt{\frac{N+\alpha-1}{4}}^2},
\]
as observed by Tzirakis \cite[Proposition 1]{Tzi2015}. Observe that, together with \eqref{eq:diamagnetic}, this gives
\begin{equation}\label{eq:trace-hardy}
  \begin{aligned}
    H_{N,\alpha}\int_{\R^N}\frac{|u(x,0)|^2}{|x|^{1-\alpha}}\dd x
    &\leq\int_{\R^{N+1}_+}y^\alpha\bigl|\nabla|u|\bigr|^2\dd x\dd y\\
    &\leq\int_{\R^{N+1}_+}y^\alpha|\nabla_\Aext u|^2\dd x\dd y.
  \end{aligned}
\end{equation}
The main result of Tzirakis \cite[Theorem 1]{Tzi2015} is a sharp interpolation between the weighted Hardy inequality and the weighted Kato inequality: for every \(\beta\in[0,N+\alpha-1)\) there is a sharp and explicit constant \(H_{N,\alpha,\beta}>0\), with \(H_{N,\alpha,0}=H_{N,\alpha}\), such that
\begin{equation}\label{eq:full-nonmag-inequality}
  \begin{aligned}
    H_{N,\alpha,\beta}\int_{\R^N}\frac{|u(x,0)|^2}{|x|^{1-\alpha}}\dd x
    &+\frac{\beta^2}{4}\int_{\R^{N+1}_+}
    \frac{y^\alpha|u|^2}{|x|^2+y^2}\dd x\dd y\\
    &\hspace{8em}
    \leq\int_{\R^{N+1}_+}y^\alpha|\nabla u|^2\dd x\dd y.
  \end{aligned}
\end{equation}
As Tzirakis observed, a convex combination of the two endpoint inequalities gives admissible constants of this type, but not the sharp ones. In the same way, combining \eqref{sharp-hardy-magnetic-inequality} with \eqref{eq:trace-hardy} does not determine the optimal intermediate trace coefficient in the magnetic setting.

This is the question we address in this paper: given \(\beta_1,\beta_2\geq0\), what is the optimal trace coefficient in the magnetic counterpart of \eqref{eq:full-nonmag-inequality}, that is, the largest constant \(K\) such that
\begin{equation}\label{eq:full-inequality}
  \begin{aligned}
    &K\int_{\R^N}\frac{|u(x,0)|^2}{|x|^{1-\alpha}}\dd x
    +\frac{\beta_1^2}{4}\int_{\R^{N+1}_+}
    \frac{y^\alpha|u|^2}{|x|^2+y^2}\dd x\dd y\\
    &\hspace{10em}+\beta_2\,\lam\int_{\R^{N+1}_+}\frac{y^\alpha|u|^2}{|x|^2}\dd x\dd y
    \leq\int_{\R^{N+1}_+}y^\alpha|\nabla_\Aext u|^2\dd x\dd y
  \end{aligned}
\end{equation}
holds for every \(u\) in the class
\[
  \mathcal C_\Aext:=C_c^\infty\pt{\overline{\R^{N+1}_+}\setminus
  \pt{\Gamma_\Sigma\times[0,\infty)}}
\]
of smooth functions that vanish near the singular set \(\Gamma_\Sigma\times[0,\infty)\) of \(\Aext\)? For \(u\in\mathcal C_\Aext\) all the integrals in \eqref{eq:full-inequality} are finite, and in terms of the Rayleigh quotient
\begin{equation}\label{eq:Rfull}
  \Rr[u]:=
  \frac{\displaystyle
    \int_{\R^{N+1}_+}y^\alpha
    \left(|\nabla_\Aext u|^2-
      \frac{\beta_1^2}{4}\frac{|u|^2}{|x|^2+y^2}-
  \beta_2\,\lam\,\frac{|u|^2}{|x|^2}\right)\dd x\dd y}
  {\displaystyle\int_{\R^N}|x|^{\alpha-1}|u(x,0)|^2\dd x},
\end{equation}
the best constant in \eqref{eq:full-inequality} is
\begin{equation}\label{eq:Kstar-def}
  \Kstar:=\inf\{\Rr[u]:\, u\in\mathcal C_\Aext,\ u(\cdot,0)\not\equiv0\}.
\end{equation}
Our main result gives an explicit formula for \(\Kstar\).

\Needspace{19\baselineskip}
\begin{samepage}
  \begin{theorem}[Sharp weighted magnetic Hardy--trace inequality]\label{thm:main}
    Let \(N\geq2\) and \(\alpha\in(-1,1)\). Let \(\Sigma\subset\Sph\) be a closed set of zero measure, let \(\a\) be a tangential vector field that is smooth on \(\Sph\setminus\Sigma\), and let \(\A(x)=\abs{x}^{-1}\a\pt{x/\abs{x}}\). Assume that
    \begin{equation}\tag{H}\label{hyp:lambda}
      N\geq3\quad\text{or}\quad\lam>0,
    \end{equation}
    let \(\beta_2\geq0\) be such that
    \[
      \mu:=\sqrt{(N-2)^2+4(1-\beta_2)\lam}>0,
    \]
    and let \(\beta_1\in[0,1+\alpha+\mu)\). Set
    \begin{equation}\label{eq:Kstar-value}
      \mathcal K(\mu)=\mathcal K(\mu,\alpha,\beta_1)
      :=
      (1-\alpha)
      \frac{
        \Gamma\pt{\frac{\alpha+1}{2}}
        \Gamma\pt{\frac{3-\alpha+\mu-\beta_1}{4}}
        \Gamma\pt{\frac{3-\alpha+\mu+\beta_1}{4}}
      }{
        \Gamma\pt{\frac{3-\alpha}{2}}
        \Gamma\pt{\frac{\alpha+1+\mu-\beta_1}{4}}
        \Gamma\pt{\frac{\alpha+1+\mu+\beta_1}{4}}
      }.
    \end{equation}
    Then,
    \begin{enumerate}[label=(\roman*)]
      \item\label{thm:main-lower} Inequality \eqref{eq:full-inequality} holds with \(K=\mathcal K(\mu)\) for every \(u\in\mathcal C_\Aext\), and all its integrals are finite. That is, \(\Kstar\geq\mathcal K(\mu)\).
      \item\label{thm:main-sharp} The constant is sharp: \(\Kstar=\mathcal K(\mu)\).
    \end{enumerate}
  \end{theorem}
\end{samepage}

When \(\lam=0\), which by \eqref{hyp:lambda} forces \(N\geq3\), one has \(\mu=N-2\); then \(\mathcal K(N-2,\alpha,\beta_1)\) is the constant \(H_{N,\alpha,\beta_1}\) of \eqref{eq:full-nonmag-inequality}, and the range of \(\beta_1\) is the same (see \cref{sec:non-magnetic}). When \(\lam>0\) and \(\beta_2<1\), one has \(\mu>N-2\), so the range \([0,1+\alpha+\mu)\) of \(\beta_1\) strictly contains the nonmagnetic range \([0,N+\alpha-1)\), and \(\mathcal K(\mu,\alpha,\beta_1)>H_{N,\alpha,\beta_1}\) for every \(\beta_1\in[0,N+\alpha-1)\) by \cref{lem:K-increasing}.

Heuristically, the constant \(\mathcal K(\mu,\alpha,\beta_1)\) is \(H_{N,\alpha,\beta_1}\) with \(N\) replaced by the ``virtual dimension'' \(D=2+\mu\). This is suggested by functions of the form \(u(x,y)=h(\abs{x},y)\psi(x/\abs{x})\), where \(\psi\) is a minimizer in \eqref{eq:lambda1}, which exists by \cref{lem:compact-embedding}: the change of unknown \(h=\abs{x}^{(D-N)/2}v\) turns \eqref{eq:full-inequality} into the nonmagnetic inequality \eqref{eq:full-nonmag-inequality} for functions of \((\abs{x},y)\), with \(\beta=\beta_1\) and with \(N\) replaced by \(D\). This identification is only formal: the nonmagnetic inequality \eqref{eq:full-nonmag-inequality} is proved in integer dimension, and a general function in \(\mathcal C_\Aext\) is not of this product form, since it involves every eigenfunction of \(-\Delta_\a\). \Cref{thm:main} shows that the formal value is nevertheless the sharp constant: the lower bound holds for every \(u\in\mathcal C_\Aext\), through a ground state representation of the deficit (see \cref{thm:deficit} below) that also shows that the constant is not attained, and the value remains sharp for smooth functions supported away from the singular cone.

\begin{remark}\label{rem:range}
  The range of \(\beta_1\) in \cref{thm:main} is exactly the condition that makes all the Gamma functions in \eqref{eq:Kstar-value} finite and positive, and \(\mathcal K(\mu)\to0\) as \(\beta_1\uparrow1+\alpha+\mu\). If \(\lam>0\), the condition \(\mu>0\) means \(\beta_2<1+(N-2)^2/(4\lam)\). The endpoint \(\beta_2=1+(N-2)^2/(4\lam)\), the critical value \(\beta_1=1+\alpha+\mu\), and the planar case \(N=2\), \(\lam=0\), which is excluded by \eqref{hyp:lambda}, are treated in the companion paper \cite{CZB}, together with the complete region of admissible interior coefficients.
\end{remark}

The proof of \cref{thm:main} is based on a ground state representation of the deficit in \eqref{eq:full-inequality} with \(K=\mathcal K(\mu)\), that is, of
\[
  \begin{aligned}
    \mathfrak D[u]:={}&\int_{\R^{N+1}_+}y^\alpha\left(
      |\nabla_\Aext u|^2-\frac{\beta_1^2}{4}\frac{|u|^2}{|x|^2+y^2}
    -\beta_2\,\lam\,\frac{|u|^2}{|x|^2}\right)\dd x\dd y\\
    &\hspace{12em}-\mathcal K(\mu)\int_{\R^N}\frac{|u(x,0)|^2}{|x|^{1-\alpha}}\dd x.
  \end{aligned}
\]

\begin{theorem}[Deficit identity and nonattainment]\label{thm:deficit}
  Under the hypotheses of \cref{thm:main}, there is a positive function \(w\) on \(\R^{N+1}_+\), depending only on \(\abs{x}\) and \(y\), smooth in \(\R^{N+1}_+\setminus(\{0\}\times(0,\infty))\) and continuous up to \(\{y=0\}\) away from \(x=0\), with \(w(x,0)=\abs{x}^{-(N+\alpha-1)/2}\), such that for every \(u\in\mathcal C_\Aext\)
  \begin{equation}\label{eq:deficit-intro}
    \begin{aligned}
      \mathfrak D[u]={}&\int_{\R^{N+1}_+}y^\alpha\,w^2\pt{\abs{\frac{x}{\abs{x}}\cdot\nabla_x\pt{\frac uw}}^2+\abs{\partial_y\pt{\frac uw}}^2}\dd x\dd y\\
      &+\int_{\R^{N+1}_+}y^\alpha\pt{\abs{\nabla_\A u-\frac{x}{\abs{x}}\pt{\frac{x}{\abs{x}}\cdot\nabla_x u}}^2-\lam\frac{\abs{u}^2}{\abs{x}^2}}\dd x\dd y,
    \end{aligned}
  \end{equation}
  and both terms on the right-hand side are nonnegative. The identity \eqref{eq:deficit-intro} extends to the completion of \(\mathcal C_\Aext\) with respect to the norm whose square is the right-hand side of \eqref{eq:full-inequality}, and on this completion \(\mathfrak D[u]=0\) only for \(u=0\). In particular, equality in \eqref{eq:full-inequality} with \(K=\mathcal K(\mu)\) is never attained.
\end{theorem}

Observe that the function \(w\) plays the role of a ground state in the following sense: if \(\psi\) is a minimizer in \eqref{eq:lambda1}, then \(w(x,y)\psi(x/\abs{x})\) solves the Euler--Lagrange equation of \(\Rr\) with critical value \(\mathcal K(\mu)\). However, this solution does not belong to the energy space, since its trace integral is infinite.

For \(\beta_1=\beta_2=0\), let \(K_{N,\a,\alpha}:=K^*_{N,\a,\alpha,0,0}\) be the best constant in
\begin{equation}\label{eq:trace-inequality}
  K_{N,\a,\alpha}
  \int_{\R^N}\frac{\abs{u(x,0)}^2}{\abs{x}^{1-\alpha}} \dd x
  \leq
  \int_{\R^{N+1}_+} y^\alpha \abs{\nabla_{\Aext}u}^2 \dd x\dd y,
  \qquad u\in\mathcal C_\Aext.
\end{equation}
Then \(K_{N,\a,\alpha}\geq H_{N,\alpha}\), and \cref{thm:main} gives \(K_{N,\a,\alpha}>H_{N,\alpha}\) whenever \(\lam>0\) (see \cref{prop:comparison}). A positive first eigenvalue may have different origins. It may come from the magnetic potential, as for \eqref{eq:AB-potential} with nonintegral flux, but also from the excluded set alone: for \(N=2\), \(\a=0\) and \(\Sigma=\{p\}\) a single point of the circle, one has \(\lam=1/4\) (see \cref{sec:example-2D}). The comparison with the nonmagnetic problem on the same angular domain, which isolates the contribution of the potential, is discussed in \cref{sec:non-magnetic}.

The link between \eqref{eq:trace-inequality} and fractional powers of the magnetic Laplacian is the extension problem introduced by Caffarelli and Silvestre \cite{CS2007} for the fractional Laplacian and extended by Stinga and Torrea \cite{StTo2010} to general nonnegative self-adjoint operators. Let \((-\Delta_\A)\) denote, as for \((-\Delta_\a)\), the Friedrichs extension of the restriction of the magnetic Laplacian \(-(\nabla-\ii\A)^2\) to \(C_c^\infty(\R^N\setminus\Gamma_\Sigma)\), that is, the operator associated with the closure of the quadratic form \(f\mapsto\int_{\R^N}|\nabla_\A f|^2\dd x\) on this space, and define \((-\Delta_\A)^\sigma\), \(\sigma\in(0,1)\), by the spectral theorem. A fixed boundary value \(f\) has many extensions to the half-space; for \(\alpha=1-2\sigma\), the least energy among them is \(\kappa_\sigma\|(-\Delta_\A)^{\sigma/2}f\|_2^2\), where \(\|\cdot\|_2\) is the norm of \(L^2(\R^N)\) and
\[
  \kappa_\sigma:=\frac{2^{1-2\sigma}\Gamma(1-\sigma)}{\Gamma(\sigma)}
\]
is the constant in the extension identity of \cite{StTo2010}.

\Needspace{22\baselineskip}
\begin{samepage}
  \begin{theorem}[Sharp fractional magnetic Hardy inequality]\label{thm:fractional}
    Let \(N\geq2\) and \(\sigma\in(0,1)\), let \(\Sigma\), \(\a\) and \(\A\) be as in \cref{thm:main}, and assume \eqref{hyp:lambda}. Set \(\alpha:=1-2\sigma\), and let
    \[
      \mathcal S_{N,\sigma,\a}
      :=
      \inf\left\{
        \frac{\displaystyle\int_{\R^N}|(-\Delta_\A)^{\sigma/2}f|^2\dd x}{\displaystyle\int_{\R^N}\abs{f}^2\abs{x}^{-2\sigma}\dd x}
        :\,
        f\in C_c^\infty\pt{\R^N\setminus\Gamma_\Sigma}\setminus\{0\}
      \right\}
    \]
    be the best constant in the fractional magnetic Hardy inequality. Then
    \begin{equation}\label{eq:S-equals-K}
      \mathcal S_{N,\sigma,\a}=\frac{K_{N,\a,\alpha}}{\kappa_\sigma}.
    \end{equation}
    Consequently, with \(\mu=\sqrt{(N-2)^2+4\lam}\), every \(f\in D(({-\Delta_\A})^{\sigma/2})\) satisfies
    \begin{equation}\label{eq:fractional-inequality}
      2^{2\sigma}
      \frac{\Gamma\pt{\frac{\mu+2+2\sigma}{4}}^2}{\Gamma\pt{\frac{\mu+2-2\sigma}{4}}^2}
      \int_{\R^N}\frac{\abs{f(x)}^2}{\abs{x}^{2\sigma}} \dd x
      \leq
      \|(-\Delta_\A)^{\sigma/2}f\|_2^2.
    \end{equation}
    The constant in \eqref{eq:fractional-inequality} equals \(\mathcal S_{N,\sigma,\a}\), so it is sharp on \(C_c^\infty(\R^N\setminus\Gamma_\Sigma)\) and on \(D((-\Delta_\A)^{\sigma/2})\). Equality in \eqref{eq:fractional-inequality} holds only for \(f=0\). For \(f\in D(({-\Delta_\A})^\sigma)\), the energy on the right also equals \(\langle f,(-\Delta_\A)^\sigma f\rangle\).
  \end{theorem}
\end{samepage}

\begin{remark}
  For \(\a=0\), \(\Sigma=\emptyset\) and \(N\geq3\), so that \(\lam=0\) and \(\mu=N-2\), the constant in \eqref{eq:fractional-inequality} becomes \(2^{2\sigma}\Gamma\pt{\frac{N+2\sigma}{4}}^2/\Gamma\pt{\frac{N-2\sigma}{4}}^2\). This is the sharp constant in the fractional Hardy inequality for \((-\Delta)^\sigma\) computed by Herbst \cite{Her1977}, which Frank, Lieb and Seiringer \cite[Proposition 4.1 and Remark 4.2]{FrLiSe2008} derived from a ground state representation that also yields remainder terms. Since \(\mu\mapsto\Gamma\pt{\frac{\mu+2+2\sigma}{4}}/\Gamma\pt{\frac{\mu+2-2\sigma}{4}}\) is increasing (see \cref{lem:K-increasing}), \eqref{eq:fractional-inequality} strictly improves this constant whenever \(\lam>0\); as for \eqref{eq:trace-inequality}, this compares with the unrestricted Euclidean problem.
\end{remark}

\medskip
\noindent\textbf{Related work.}
Beyond the homogeneous setting, Cazacu and Krej\v{c}i\v{r}\'ik \cite{CK2016} established magnetic Hardy inequalities in every dimension and used them to describe the large-time decay of the heat semigroup, which for compactly supported fields is governed by an eigenvalue problem for a magnetic operator on the sphere. Fanelli, Krej\v{c}i\v{r}\'ik, Laptev and Vega \cite[Theorem 1]{FKLV2020} extended the Laptev--Weidl inequality \eqref{eq:hardy-magnetic2} to the Aharonov--Bohm potential in \(\R^N\), \(N\geq3\), which is singular along the subspace \(\{x_1=x_2=0\}\): they proved, with sharp constants,
\[
  \pt{\frac{N-2}{2}}^2\int_{\R^N}\frac{\abs{f(x)}^2}{\abs{x}^2}\dd x
  +\dist(\bar\a,\Z)^2\int_{\R^N}\frac{\abs{f(x)}^2}{x_1^2+x_2^2}\dd x
  \leq
  \int_{\R^N}\abs{\nabla_\A f(x)}^2\dd x.
\]
For \(N=3\), this potential is the one of \cref{sec:example-AB3D}, for which \(\lam=b(b+1)\) with \(b=\dist(\bar\a,\Z)\), so \eqref{eq:hardy-magnetic} has the coefficient \(\frac14+b(b+1)\) in front of \(\abs{x}^{-2}\). The two inequalities are not comparable: the weight \((x_1^2+x_2^2)^{-1}\) is larger than \(\abs{x}^{-2}\), but its coefficient \(b^2\) is smaller than \(b(b+1)\). Fanelli and Kova\v{r}\'ik \cite{FK2024} proved quantitative Hardy inequalities for two-dimensional magnetic Hamiltonians, and the recent preprint of Kova\v{r}\'ik and Rossaro \cite{KR2026} studies magnetic Hardy inequalities with singular integral weights. In nonlinear magnetic interpolation problems, attainment and symmetry are central questions: Dolbeault, Esteban, Laptev and Loss \cite{DELL2018} obtained interpolation inequalities and spectral estimates for magnetic operators, Bonheure, Dolbeault, Esteban, Laptev and Loss \cite{BDELL2020,BDELL2021} studied symmetry and symmetry breaking for Aharonov--Bohm fields, and Abatangelo and Terracini \cite{AT2011} studied nonlinear Schr\"odinger equations with singular electromagnetic potentials and critical exponent. Our problem is quadratic, and the situation is different: the first eigenfunction of \(-\Delta_\a\) exists by compactness, but \eqref{eq:full-inequality} has no extremal functions, and the sharp constant is only approached by sequences of functions in \(\mathcal C_\Aext\) (see \cref{prop:Kstar-le-Kred}).

For the nonmagnetic inequalities \eqref{eq:full-nonmag-inequality}, Alvino, Volpicelli and Ferone \cite{AVF2012} identified, for \(\alpha=0\), the optimal trace remainder below the critical interior coefficient, Tzirakis \cite{Tzi2022} obtained series expansions for weighted Finsler versions of these inequalities, and Nguyen \cite[Theorem 1.1]{Ngu2016} proved the analogous sharp inequalities on polyhedral convex cones, with \(y^\alpha\) replaced by a power of the distance to the boundary of the cone. For fractional energies, diamagnetic bounds give the Euclidean constant as a baseline when the potential is locally square-integrable \cite[(2.3)--(2.4)]{FrLiSe2008}; this excludes the planar Aharonov--Bohm potential and the three-dimensional examples of \cref{sec:examples}. Expanding in eigenfunctions of \((-\Delta_\a)\) decomposes \((-\Delta_\A)\) into Bessel operators in the radial variable \cite{CZB}; the one associated with the first eigenvalue has dimension parameter \(2+\mu\), and its sharp fractional Hardy constant is the Gamma quotient in \eqref{eq:fractional-inequality}. For a single Bessel operator this constant is a nonmagnetic fact: it was computed by Yafaev \cite[(2.4) and (2.26)]{Yaf1999} for the dimension parameters arising from spherical harmonics, by a one-dimensional argument \cite[Lemmas 2.2, 2.6 and 2.7]{Yaf1999} that applies to every real value of this parameter, and it also follows from the ground state representation of Bogdan and Merz \cite[Theorem 3.2]{BM2024} for fractional powers of Bessel operators, which holds for real values of the dimension parameter. The same authors obtained heat kernel bounds for the fractional Laplacian with Hardy potential in angular momentum channels \cite{BM2026}. Our contribution on the fractional side is to prove the extension principle on the form domain of \((-\Delta_\A)^\sigma\), to prove sharpness on \(C_c^\infty(\R^N\setminus\Gamma_\Sigma)\), and to connect the fractional inequality with the whole interior--trace family. This spectral setting should also be distinguished from phase-modified Gagliardo models, such as the one in \cite{LJG2019}, which are not identified here with powers of \((-\Delta_\A)\).

The paper is organized as follows. In \Cref{sec:preliminaries} we collect a few technical lemmata. \Cref{sec:reducing-PDE} reduces the Euler--Lagrange equation of \(\Rr\) to an ODE. \Cref{sec:general-ode} constructs an explicit positive solution of this ODE, which determines the sharp constant, and \cref{sec:main-proof} proves \cref{thm:main,thm:deficit}. \Cref{sec:non-magnetic} compares \(\mathcal K(\mu)\) with the constant of the nonmagnetic problem and characterizes \(\lam=0\) in terms of capacity and gauges (\cref{prop:lambda1-zero}). \Cref{sec:fractional} proves the extension principle and \cref{thm:fractional}. \Cref{sec:examples} illustrates the hypotheses of \cref{thm:main} with magnetic examples: in the plane we revisit the Aharonov--Bohm potential \eqref{eq:AB-potential}, and in three dimensions we study the Aharonov--Bohm potential and a nontrivial Dirac monopole. The appendices collect the hypergeometric identities and the facts on capacity used in the proofs.

\Needspace{6\baselineskip}
\section{Preliminaries}\label{sec:preliminaries}
The compactness of the embedding of \(H_\a^1(\Sph)\) into \(L^2(\Sph)\) is standard, but we give a proof for completeness. In particular, the infimum in \eqref{eq:lambda1} is attained, so a first eigenfunction of \((-\Delta_\a)\) is available for the separation of variables in \cref{sec:magnetic-reduction}.

\begin{lemma}[Compactness of the angular form domain]\label{lem:compact-embedding}
  Under the angular assumptions of \cref{thm:main}, the magnetic energy in \eqref{eq:HA1}, defined on \(C_c^\infty(\Sph\setminus\Sigma)\), is closable in \(L^2(\Sph)\), and the embedding \(H_\a^1(\Sph)\hookrightarrow L^2(\Sph)\) is compact. Consequently, \(L^2(\Sph)\) has an orthonormal basis \(\{\psi_j\}_{j=1}^\infty\) of eigenfunctions of \((-\Delta_\a)\), with eigenvalues \(\lam=\lambda_1\leq\lambda_2\leq\cdots\to\infty\) repeated according to multiplicity, and the infimum in \eqref{eq:lambda1} is attained at \(\psi_1\).
\end{lemma}

\begin{proof}[Proof of \cref{lem:compact-embedding}]
  Put \(\Omega:=\Sph\setminus\Sigma\). Since \(\abs\Sigma=0\), \(C_c^\infty(\Omega)\) is dense in \(L^2(\Sph)\), and the form is that of \((-\Delta_\a)\) on \(C_c^\infty(\Omega)\), which is symmetric and nonnegative; hence it is closable, and \(H_\a^1(\Sph)\subset L^2(\Sph)\).

  \noindent\underline{\emph{Estimate near \(\Sigma\).}} Let \(\psi\in C_c^\infty(\Omega)\). Since \(\a\) is real, \(\nabla_{\Sph}\abs\psi=\operatorname{Re}(\overline\psi\,\nabla_\a\psi)/\abs\psi\) where \(\psi\neq0\), which gives the diamagnetic inequality \(\abs{\nabla_{\Sph}\abs\psi}\leq\abs{\nabla_\a\psi}\) almost everywhere~\cite{LiebLoss1997}. Let \(p:=2N/(N-1)>2\), which is below the Sobolev exponent of the \((N-1)\)-dimensional sphere, so that \(H^1(\Sph)\subset L^p(\Sph)\). Since \(1-2/p=1/N\), H\"older's inequality gives, for every measurable \(G\subset\Sph\),
  \begin{equation}\label{eq:small-set}
    \int_G\abs\psi^2\dd\omega
    \leq \abs G^{1/N}\norm{\abs\psi}_{L^p}^2
    \leq C\abs G^{1/N}\norm{\psi}_{H_\a^1}^2,
  \end{equation}
  and by density \eqref{eq:small-set} holds for every \(\psi\in H_\a^1(\Sph)\).

  \noindent\underline{\emph{Total boundedness.}} Let \(\mathcal B\) be the unit ball of \(H_\a^1(\Sph)\) and let \(\varepsilon>0\). If \(\Sigma=\emptyset\), then \(\a\) is bounded, \(\mathcal B\) is bounded in \(H^1(\Sph)\), and the Rellich--Kondrachov theorem shows that \(\mathcal B\) is totally bounded in \(L^2(\Sph)\). Otherwise, let \(G_t:=\{\omega\in\Sph:\dist(\omega,\Sigma)<t\}\), with \(\dist\) the geodesic distance. Since \(\Sigma\) is closed, \(\abs{G_t}\to\abs\Sigma=0\) as \(t\to0\), so we may fix \(\delta>0\) with \(C\abs{G_{2\delta}}^{1/N}<\varepsilon^2\). Let
  \[
    \phi:=\min\bigg\{1,\pt{\frac{\dist(\cdot,\Sigma)}{\delta}-1}_+\bigg\},
  \]
  which is Lipschitz, vanishes in \(G_\delta\), and equals \(1\) outside \(G_{2\delta}\). The field \(\a\) is bounded on the compact set \(\Sph\setminus G_\delta\subset\Omega\), so for \(\psi\in C_c^\infty(\Omega)\) the identity \(\nabla_{\Sph}(\phi\psi)=\phi(\nabla_\a\psi+\ii\a\psi)+\psi\nabla_{\Sph}\phi\) gives \(\norm{\phi\psi}_{H^1(\Sph)}\leq C_\delta\norm{\psi}_{H_\a^1}\). Since \(\phi\psi_k\to\phi\psi\) in \(L^2(\Sph)\) whenever \(\psi_k\to\psi\) in \(H_\a^1(\Sph)\), this bound extends by density to every \(\psi\in H_\a^1(\Sph)\). Hence \(\{\phi\psi:\psi\in\mathcal B\}\) is bounded in \(H^1(\Sph)\), and by the Rellich--Kondrachov theorem there are \(g_1,\ldots,g_n\in L^2(\Sph)\) such that every \(\psi\in\mathcal B\) satisfies \(\norm{\phi\psi-g_j}_{L^2}<\varepsilon\) for some \(j\). Since \(0\leq1-\phi\leq\chi_{G_{2\delta}}\), \eqref{eq:small-set} gives \(\norm{(1-\phi)\psi}_{L^2}<\varepsilon\), and therefore \(\norm{\psi-g_j}_{L^2}<2\varepsilon\). Thus \(\mathcal B\) is totally bounded in \(L^2(\Sph)\), and the embedding is compact.
\end{proof}

We use spherical coordinates \(x=\rho\omega\) with \(\rho=\abs{x}>0\) and \(\omega\in\Sph\), so that \(\dd x=\rho^{N-1}\dd\rho\dd\omega\). For a function \(u\) on \(\R^{N+1}_+\) and fixed \((\rho,y)\), we write \(u(\rho\,\cdot,y)\) for the function \(\omega\mapsto u(\rho\omega,y)\) on \(\Sph\); for functions on \(\Sph\), \(\|\cdot\|_2\) denotes the norm of \(L^2(\Sph)\). For functions of \((\rho,y)\in(0,\infty)^2\) we write \(\nabla_{(\rho,y)}\) and \(\diver_{(\rho,y)}\) for the gradient and the divergence in these variables. Since \(\A(x)=\rho^{-1}\a(\omega)\) and \(\a\) is tangential,
\begin{equation}\label{eq:grad-spherical}
  \nabla_x-\ii \A=\omega\partial_\rho+\frac{1}{\rho}\nabla_{\a},
  \qquad
  \abs{\nabla_{\Aext}u}^2=\abs{\partial_\rho u}^2+\abs{\partial_yu}^2+\frac{1}{\rho^2}\abs{\nabla_\a u}^2,
\end{equation}
where \(\nabla_\a\) acts on the variable \(\omega\); the second identity holds because \(\nabla_\a u\) is tangent to \(\Sph\), hence orthogonal to \(\omega\).

\begin{lemma}[Separated-variable identities]\label{lem:separated}
  Let \(h:(0,\infty)\times[0,\infty)\to\C\) be continuous and locally
  Lipschitz in \((0,\infty)^2\), and let
  \(\varphi\in C_c^\infty(\Sph\setminus\Sigma)\).
  For \(u(\rho\omega,y)=h(\rho,y)\varphi(\omega)\), the following
  identities hold in \([0,\infty]\):
  \begin{align*}
    \int_{\R^{N+1}_+}y^\alpha|\nabla_\Aext u|^2\dd x\dd y&=\|\varphi\|_2^2
    \int_0^\infty\!\int_0^\infty\rho^{N-1}y^\alpha|\nabla_{(\rho,y)}h|^2\dd\rho\dd y\\
    &\quad+\|\nabla_\a\varphi\|_2^2
    \int_0^\infty\!\int_0^\infty\rho^{N-3}y^\alpha|h|^2\dd\rho\dd y,\\
    \int_{\R^{N+1}_+}\frac{y^\alpha|u|^2}{|x|^2+y^2}\dd x\dd y&=\|\varphi\|_2^2
    \int_0^\infty\!\int_0^\infty\frac{\rho^{N-1}y^\alpha|h|^2}{\rho^2+y^2}\dd\rho\dd y,\\
    \int_{\R^{N+1}_+}\frac{y^\alpha|u|^2}{|x|^2}\dd x\dd y&=\|\varphi\|_2^2
    \int_0^\infty\!\int_0^\infty\rho^{N-3}y^\alpha|h|^2\dd\rho\dd y,\\
    \int_{\R^N}|x|^{\alpha-1}|u(x,0)|^2\dd x&=\|\varphi\|_2^2
    \int_0^\infty\rho^{N+\alpha-2}|h(\rho,0)|^2\dd\rho.
  \end{align*}
\end{lemma}

\begin{proof}[Proof of Lemma~\ref{lem:separated}]
  Since \(h\) does not depend on \(\omega\), \eqref{eq:grad-spherical} gives \(\abs{\nabla_{\Aext}u}^2=\abs{\varphi}^2\abs{\nabla_{(\rho,y)}h}^2+\rho^{-2}\abs{h}^2\abs{\nabla_\a\varphi}^2\) almost everywhere. The identities follow from \(\dd x=\rho^{N-1}\dd\rho\dd\omega\) and Tonelli's theorem.
\end{proof}

The following monotonicity property of the constant \eqref{eq:Kstar-value} is used in \cref{sec:non-magnetic,sec:fractional}.
\begin{lemma}[Monotonicity of the trace coefficient]\label{lem:K-increasing}
  Let \(\alpha\in(-1,1)\) and \(\beta_1\geq0\). The function
  \(\mu\mapsto\mathcal K(\mu,\alpha,\beta_1)\) is strictly increasing on
  \((\beta_1-1-\alpha,\infty)\).
\end{lemma}
\begin{proof}[Proof of Lemma~\ref{lem:K-increasing}]
  Apart from a positive constant independent of \(\mu\), the coefficient is
  the product, for the two signs, of
  \[
    G_\pm(\mu)=
    \frac{\Gamma((3-\alpha\pm\beta_1+\mu)/4)}
    {\Gamma((1+\alpha\pm\beta_1+\mu)/4)}.
  \]
  All arguments are positive, and each numerator argument exceeds the
  denominator argument by \((1-\alpha)/2\). The strict increase of the
  digamma function then gives \((\log G_\pm)'(\mu)>0\).
\end{proof}

\Needspace{6\baselineskip}
\section{Reduction to an ODE}\label{sec:reducing-PDE}
\Needspace{5\baselineskip}
\subsection{Spherical reduction}\label{sec:magnetic-reduction}
In this subsection we separate the angular variable and derive the equation satisfied by the radial–vertical profile of the first angular mode.

By \cref{lem:compact-embedding}, \((-\Delta_\a)\) has a first eigenfunction \(\psi_1\), that is, \(-\Delta_\a\psi_1=\lam\psi_1\); it is smooth on \(\Sph\setminus\Sigma\) by elliptic regularity. Since the magnetic Laplacian in the variable \(x\) is \(\partial_{\rho\rho}+\rho^{-1}(N-1)\partial_\rho+\rho^{-2}\Delta_\a\), a function \(u(\rho\omega,y)=w(\rho,y)\psi_1(\omega)\) formally solves the Euler--Lagrange equation of \(\Rr\), with critical value \(K\), if \(w\) solves
\begin{equation}\label{eq:mode-equation}
  -\diver_{(\rho,y)}\pt{\rho^{N-1}y^\alpha\nabla_{(\rho,y)}w}
  =
  \rho^{N-1}y^\alpha\pt{\frac{\beta_1^2}{4(\rho^2+y^2)}+\frac{(\beta_2-1)\lam}{\rho^2}}w
  \quad\text{ in }\;(0,\infty)^2,
\end{equation}
together with
\begin{equation}\label{eq:boundary-mode}
  -\lim_{y\to0^+}\rho^{N-1}y^\alpha\partial_yw(\rho,y)
  =
  K\rho^{N+\alpha-2}w(\rho,0),
  \qquad\rho>0.
\end{equation}

\subsection{Self-similar reduction}\label{sec:selfsimilar-reduction}
In this subsection we reduce \eqref{eq:mode-equation}--\eqref{eq:boundary-mode} to an ordinary differential equation in the variable \(t=y/\rho\). Following \cite{Tzi2015}, write
\begin{equation}\label{eq:w-selfsimilar}
  \begin{aligned}
    w(\rho,y) & =(\rho^2+y^2)^{\gamma/2}\rho^\delta V(y/\rho)
    =\rho^s\Phi(y/\rho),\quad
    s=\gamma+\delta,\\
    \Phi(t) &=(1+t^2)^{\gamma/2}V(t).
  \end{aligned}
\end{equation}
Substitution in \eqref{eq:mode-equation} gives
\begin{equation}\label{eq:Phi-ODE}
  \begin{aligned}
    & (1+t^2)\Phi''+
    \left(\frac\alpha t+(3-N-2s)t\right)\Phi'\\
    &\hspace{8em}+\left(s(s+N-2)+(\beta_2-1)\lam+
    \frac{\beta_1^2}{4(1+t^2)}\right)\Phi=0.
  \end{aligned}
\end{equation}
The homogeneity of the trace quotient suggests the choice
\begin{equation}\label{eq:special-choice}
  \gamma=-\frac{\beta_1}{2},\qquad
  \delta=-\frac{(N+\alpha-1-\beta_1)}{2},\qquad
  s=-\frac{(N+\alpha-1)}{2}.
\end{equation}
Substituting \(\Phi=(1+t^2)^{\gamma/2}V\) in \eqref{eq:Phi-ODE} and dividing by \((1+t^2)^{\gamma/2}\), we obtain a second-order equation for \(V\) whose zeroth-order coefficient is
\[
  \begin{gathered}
    C_0(t):=s(s+N-2)+(\beta_2-1)\lam+\frac{C_1}{1+t^2}+\frac{C_2t^2}{1+t^2},\\
    C_1:=\gamma(1+\alpha)+\frac{\beta_1^2}{4},
    \qquad
    C_2:=\gamma(\gamma+2-N-2s).
  \end{gathered}
\]
The choice \eqref{eq:special-choice} gives \(C_1=C_2=\beta_1^2/4-\beta_1(1+\alpha)/2\), so \(C_0(t)=s(s+N-2)+(\beta_2-1)\lam+C_1\) does not depend on \(t\), and in terms of the parameters can be written as
\begin{equation}\label{eq:C0-special}
  C_0(t)=C_0:=-\frac14(N+\alpha-1-\beta_1)(N-\alpha-3+\beta_1)
  -(1-\beta_2)\lam.
\end{equation}
For the boundary condition we use
\begin{equation}\label{eq:Phi-prime}
  \Phi'(t)=(1+t^2)^{\gamma/2}
  \left(V'(t)+\frac{\gamma t}{1+t^2}V(t)\right).
\end{equation}
The common factor \(\rho^{s+N+\alpha-2}\) in the weighted normal
derivative and the boundary value cancels. Thus the reduced problem is
\begin{equation}\label{eq:V-ODE-first}
  \left\{
    \begin{aligned}
      &(1+t^2)V''+\left(\frac\alpha t+(2+\alpha-\beta_1)t\right)V'+C_0V=0,
      \quad t>0,\\
      &-\lim_{t\to0^+}t^\alpha
      \left(V'(t)+\frac{\gamma t}{1+t^2}V(t)\right)=KV(0).
    \end{aligned}\right.
  \end{equation}

  We now compute the integrals of \cref{lem:separated} for functions \(\eta(\rho)w(\rho,y)\), with \(w\) as in \eqref{eq:w-selfsimilar} and a radial cutoff \(\eta\). In \cref{sec:upper-bound}, these functions are used in the proof of the upper bound \(\Kstar\leq\mathcal K(\mu)\).

  For a function \(V:(0,\infty)\to\R\) with a limit \(V(0^+)\), we let
  \[
    P_1[V]
    :=
    \int_0^\infty t^\alpha(1+t^2)^\gamma V^2\dd t,
    \qquad
    T[V]
    :=
    V(0^+)^2,
  \]
  where \(\gamma\) is given by \eqref{eq:special-choice}.
  \begin{lemma}[Radial factorization identities]\label{lem:radial-factor}
    Let \(w\) be given by \eqref{eq:w-selfsimilar} and \eqref{eq:special-choice}, with \(V\) continuous on \([0,\infty)\), and let \(\eta:(0,\infty)\to[0,\infty)\) be measurable. Then, as identities in \([0,\infty]\),
    \begin{align*}
      \int_0^\infty\int_0^\infty\rho^{N-3}y^\alpha\eta(\rho)^2w^2\dd\rho\dd y
      &=
      P_1[V]\int_0^\infty\eta^2\frac{\dd\rho}{\rho},
      \\
      \int_0^\infty\int_0^\infty\rho^{N-1}y^\alpha\eta(\rho)^2w^2\dd\rho\dd y
      &=
      P_1[V]\int_0^\infty\eta^2\rho\dd\rho,
      \\
      \int_0^\infty\rho^{N+\alpha-2}\eta(\rho)^2w(\rho,0)^2\dd\rho
      &=
      T[V]\int_0^\infty\eta^2\frac{\dd\rho}{\rho}.
    \end{align*}
  \end{lemma}

  \begin{proof}[Proof of Lemma~\ref{lem:radial-factor}]
    We have \(w^2=\rho^{2s}(1+t^2)^\gamma V^2\) with \(t=y/\rho\). In the first two integrals we substitute \(y=\rho t\), so that \(y^\alpha\dd y=\rho^{\alpha+1}t^\alpha\dd t\), and use Tonelli's theorem. Since \(2s=-(N+\alpha-1)\), the resulting powers of \(\rho\) are \(\rho^{N-3+\alpha+1+2s}=\rho^{-1}\) and \(\rho^{N-1+\alpha+1+2s}=\rho\). In the last one, \(w(\rho,0)=\rho^sV(0)\) and \(\rho^{N+\alpha-2+2s}=\rho^{-1}\).
  \end{proof}

  \section{The admissible self-similar profile}\label{sec:general-ode}
  In this section we construct an explicit positive solution of \eqref{eq:V-ODE-first} with \(P_1[V]<\infty\) and \(T[V]=1\), and compute the corresponding value of \(K\).

  We retain the hypotheses of \cref{thm:main}, in particular \(\mu>0\) and \(0\leq\beta_1<1+\alpha+\mu\).
  The change of variable \(z=-t^2\) puts \eqref{eq:V-ODE-first} into the
  hypergeometric equation \eqref{eq:gauss-ode}, with parameters
  \begin{equation}\label{eq:ABC}
    A_\pm=\frac{(1+\alpha-\beta_1\pm\mu)}{4},\qquad
    C=\frac{(1+\alpha)}{2},\qquad q_*=1+\frac\mu2.
  \end{equation}
  Indeed,
  \begin{equation}\label{eq:C0-A}
    C_0=4A_+A_-=\frac{(1+\alpha-\beta_1)^2-\mu^2}{4}.
  \end{equation}
  Since \(\beta_1<1+\alpha+\mu\) and \(-1<\alpha<1\),
  \begin{equation}\label{eq:ABC-positive}
    \begin{gathered}
      A_+>0,\qquad C-A_-=\frac{1+\alpha+\beta_1+\mu}{4}>0,\\
      q_*-A_+-(C-A_-)=1-C=\frac{1-\alpha}{2}>0.
    \end{gathered}
  \end{equation}
  Define
  \begin{equation}\label{eq:C-fast}
    C_{\mathrm{fast}}
    =\frac{\Gamma(1-A_-)\Gamma(A_+-C+1)}
    {\Gamma(q_*)\Gamma(1-C)}>0
  \end{equation}
  and
  \begin{equation}\label{eq:Vadm}
    \Vadm(t)=C_{\mathrm{fast}}(1+t^2)^{-A_+}
    \F\left(A_+,C-A_-;q_*;\frac{1}{1+t^2}\right),\qquad t\geq0.
  \end{equation}
  The identities used below are recorded in \cref{sec:hypergeometric}.

  \Needspace{11\baselineskip}
  \begin{proposition}[Admissible profile and boundary flux]\label{prop:explicit-Kred}
    The function \(\Vadm\) solves \eqref{eq:V-ODE-first} and has the following properties.
    \begin{enumerate}[label=(\roman*)]
      \item\label{it:Vadm-asy} As \(t\to\infty\),
        \[
          \Vadm(t)=C_{\mathrm{fast}}t^{-2A_+}(1+O(t^{-2})),\qquad
          \Vadm'(t)=-2A_+C_{\mathrm{fast}}t^{-2A_+-1}(1+O(t^{-2})).
        \]
      \item\label{it:Vadm-local} As \(t\to0^+\),
        \[
          \Vadm(t)=1+\kappa_{\mathrm{adm}}t^{1-\alpha}+O(t^2),\qquad
          \Vadm'(t)=(1-\alpha)\kappa_{\mathrm{adm}}t^{-\alpha}+O(t),
        \]
        where
        \begin{equation}\label{eq:kappa-adm-expanded}
          \kappa_{\mathrm{adm}}
          =-\frac{\Gamma(C)\Gamma(A_+-C+1)\Gamma(1-A_-)}
          {\Gamma(2-C)\Gamma(A_+)\Gamma(C-A_-)}.
        \end{equation}
      \item\label{it:Vadm-finite} \(P_1[\Vadm]<\infty\) and \(T[\Vadm]=1\).
      \item\label{it:Vadm-pos} \(\Vadm(t)>0\) for every \(t\geq0\).
      \item\label{it:Vadm-K} The boundary coefficient is
        \[
          K(\Vadm):=-\lim_{t\to0^+}t^\alpha
          \left(\Vadm'(t)+\frac{\gamma t}{1+t^2}\Vadm(t)\right)
          =-(1-\alpha)\kappa_{\mathrm{adm}}=\mathcal K(\mu,\alpha,\beta_1).
        \]
    \end{enumerate}
  \end{proposition}

  \begin{proof}[Proof of Proposition~\ref{prop:explicit-Kred}]
    Kummer's solution \eqref{eq:kummer-infinity} with \((a,b,c)=(A_+,A_-,C)\) and \(z=-t^2\), where \(A_+-A_-+1=q_*\), shows that \(t^{-2A_+}\F(A_+,A_+-C+1;q_*;-t^{-2})\) solves the differential equation in \eqref{eq:V-ODE-first}. Pfaff's identity~\eqref{eq:pfaff} rewrites it as \eqref{eq:Vadm} divided by \(C_{\mathrm{fast}}\), so \(\Vadm\) solves the same equation. At \(t=0\) the argument of \(\F\) in \eqref{eq:Vadm} equals \(1\), and Gauss's formula \eqref{eq:gauss-sum}, with \(c-a-b=1-C>0\), gives
    \[
      \Vadm(0)=C_{\mathrm{fast}}\frac{\Gamma(q_*)\Gamma(1-C)}{\Gamma(1-A_-)\Gamma(A_+-C+1)}=1;
    \]
    this is the reason for the choice \eqref{eq:C-fast}.

    \noindent\underline{\emph{Proof of \ref{it:Vadm-asy}.}} Write \(\Vadm(t)=C_{\mathrm{fast}}(1+t^2)^{-A_+}g\pt{(1+t^2)^{-1}}\) with \(g(x):=\F(A_+,C-A_-;q_*;x)\), which is analytic near \(x=0\) with \(g(0)=1\). As \(t\to\infty\), \((1+t^2)^{-1}=O(t^{-2})\), so \(g\pt{(1+t^2)^{-1}}=1+O(t^{-2})\), and \((1+t^2)^{-A_+}=t^{-2A_+}(1+O(t^{-2}))\). This gives the expansion of \(\Vadm\). Differentiating,
    \[
      \Vadm'(t)
      =-\frac{2t\,C_{\mathrm{fast}}}{(1+t^2)^{A_++1}}\pt{A_+\,g\pt{\tfrac1{1+t^2}}+\frac{g'\pt{\frac1{1+t^2}}}{1+t^2}}
      =-2A_+C_{\mathrm{fast}}t^{-2A_+-1}(1+O(t^{-2})).
    \]

    \noindent\underline{\emph{Proof of \ref{it:Vadm-local}.}} We apply \eqref{eq:connection-one} with
    \[
      a=A_+,\qquad b=C-A_-,\qquad c=q_*,\qquad z=\frac1{1+t^2},\qquad 1-z=\frac{t^2}{1+t^2},
    \]
    which is valid for \(t>0\) because \(c-a-b=1-C\in(0,1)\) is not an integer. Since \(A_+-A_-=\mu/2\), we have
    \[
      c-a=1-A_-,\quad c-b=A_+-C+1,\quad a+b-c+1=C,\quad c-a-b+1=2-C.
    \]
    Moreover \((1-z)^{1-C}=t^{1-\alpha}(1+t^2)^{-(1-\alpha)/2}\), because \(2(1-C)=1-\alpha\). Multiplying \eqref{eq:connection-one} by \(C_{\mathrm{fast}}(1+t^2)^{-A_+}\) therefore gives
    \begin{equation}\label{eq:Vadm-two-branches}
      \Vadm(t)=c_1G_1(t)+c_2\,t^{1-\alpha}G_2(t),\qquad t>0,
    \end{equation}
    where
    \[
      \begin{aligned}
        G_1(t)&:=(1+t^2)^{-A_+}\F\pt{A_+,C-A_-;C;\tfrac{t^2}{1+t^2}},\\
        G_2(t)&:=(1+t^2)^{-A_+-(1-\alpha)/2}\F\pt{1-A_-,A_+-C+1;2-C;\tfrac{t^2}{1+t^2}},
      \end{aligned}
    \]
    and
    \[
      c_1=C_{\mathrm{fast}}\frac{\Gamma(q_*)\Gamma(1-C)}{\Gamma(1-A_-)\Gamma(A_+-C+1)},
      \qquad
      c_2=C_{\mathrm{fast}}\frac{\Gamma(q_*)\Gamma(C-1)}{\Gamma(A_+)\Gamma(C-A_-)}.
    \]
    By the definition \eqref{eq:C-fast} of \(C_{\mathrm{fast}}\), \(c_1=1\). Since \(\Gamma(C)=(C-1)\Gamma(C-1)\) and \(\Gamma(2-C)=(1-C)\Gamma(1-C)\), we have \(\Gamma(C-1)/\Gamma(1-C)=-\Gamma(C)/\Gamma(2-C)\), and therefore
    \[
      c_2=\frac{\Gamma(1-A_-)\Gamma(A_+-C+1)\Gamma(C-1)}{\Gamma(1-C)\Gamma(A_+)\Gamma(C-A_-)}
      =-\frac{\Gamma(C)\Gamma(A_+-C+1)\Gamma(1-A_-)}{\Gamma(2-C)\Gamma(A_+)\Gamma(C-A_-)}
      =\kappa_{\mathrm{adm}}.
    \]
    Since \(C\) and \(2-C\) are not integers, \(G_1\) and \(G_2\) are analytic functions of \(t^2\) near \(t=0\) with \(G_1(0)=G_2(0)=1\). Hence \(G_i(t)=1+O(t^2)\) and \(G_i'(t)=O(t)\) for \(i=1,2\). Inserting this in \eqref{eq:Vadm-two-branches},
    \[
      \begin{aligned}
        \Vadm(t)&=1+\kappa_{\mathrm{adm}}t^{1-\alpha}+O(t^2)+O(t^{3-\alpha}),\\
        \Vadm'(t)&=(1-\alpha)\kappa_{\mathrm{adm}}t^{-\alpha}+O(t)+O(t^{2-\alpha}),
      \end{aligned}
    \]
    and since \(\alpha<1\), the terms \(O(t^{3-\alpha})\) and \(O(t^{2-\alpha})\) are \(O(t^2)\) and \(O(t)\), respectively, as \(t\to0^+\).

    \noindent\underline{\emph{Proof of \ref{it:Vadm-finite}.}} By \ref{it:Vadm-local}, \(\Vadm(0^+)=1\), so \(T[\Vadm]=1\). The integrand \(t^\alpha(1+t^2)^\gamma\Vadm^2\) of \(P_1[\Vadm]\) is bounded by a constant times \(t^\alpha\) near zero, which is integrable because \(\alpha>-1\). By \ref{it:Vadm-asy} and \(\gamma=-\beta_1/2\), it is of order \(t^{\alpha-\beta_1-4A_+}=t^{-1-\mu}\) at infinity, which is integrable because \(\mu>0\). Hence \(P_1[\Vadm]<\infty\).

    \noindent\underline{\emph{Proof of \ref{it:Vadm-pos}.}} For \(t\geq0\) the argument \(x=(1+t^2)^{-1}\) of \(\F\) in \eqref{eq:Vadm} lies in \((0,1]\). By \eqref{eq:ABC-positive}, every term of the series \eqref{def:hypergeometric} defining \(\F(A_+,C-A_-;q_*;x)\) is positive, and its \(n\)-th coefficient is of order \(n^{-1-(1-C)}\), so the series also converges at \(x=1\). Hence \(\F(A_+,C-A_-;q_*;x)\geq1\), and since \(C_{\mathrm{fast}}>0\), \(\Vadm(t)>0\) for every \(t\geq0\).

    \noindent\underline{\emph{Proof of \ref{it:Vadm-K}.}} By \ref{it:Vadm-local}, \(t^\alpha\Vadm'(t)=(1-\alpha)\kappa_{\mathrm{adm}}+O(t^{1+\alpha})\) and \(t^\alpha\frac{\gamma t}{1+t^2}\Vadm(t)=O(t^{1+\alpha})\) as \(t\to0^+\). Since \(\alpha>-1\), this gives \(K(\Vadm)=-(1-\alpha)\kappa_{\mathrm{adm}}\). By \eqref{eq:ABC},
    \[
      \begin{aligned}
        C&=\frac{1+\alpha}{2}, & 2-C&=\frac{3-\alpha}{2},\\
        A_+&=\frac{1+\alpha-\beta_1+\mu}{4}, & A_+-C+1&=\frac{3-\alpha-\beta_1+\mu}{4},\\
        C-A_-&=\frac{1+\alpha+\beta_1+\mu}{4}, & 1-A_-&=\frac{3-\alpha+\beta_1+\mu}{4},
      \end{aligned}
    \]
    and inserting these values in \eqref{eq:kappa-adm-expanded} gives \(-(1-\alpha)\kappa_{\mathrm{adm}}=\mathcal K(\mu,\alpha,\beta_1)\), see \eqref{eq:Kstar-value}.
  \end{proof}

  \section{Proof of Theorems \texorpdfstring{\ref{thm:main} and \ref{thm:deficit}}{1.1 and 1.3}}\label{sec:main-proof}
  Throughout this section we assume \eqref{hyp:lambda}: \(N\geq 3\) or \(\lam>0\), and the remaining hypotheses of \cref{thm:main}, namely
  \[
    \beta_2\geq0,\qquad
    \mu>0,\qquad
    0\leq\beta_1<1+\alpha+\mu.
  \]
  We let \(\Vadm\) be given by \eqref{eq:Vadm}, so that \(K(\Vadm)=\mathcal K(\mu)\) by \cref{prop:explicit-Kred}, and
  \begin{equation}\label{eq:w-first-self}
    w(\rho,y):=(\rho^2+y^2)^{\gamma/2}\rho^\delta\Vadm\pt{\frac{y}{\rho}},
  \end{equation}
  with \(\gamma\), \(\delta\) as in \eqref{eq:special-choice}. We call \(w\) the profile.

  \subsection{A ground state representation}

  \begin{lemma}[Reduced ground-state representation]\label{lem:reduced-GSR}
    The function \(w\) in~\eqref{eq:w-first-self} is positive on \((0,\infty)\times[0,\infty)\). Let \(h:(0,\infty)\times[0,\infty)\to\C\) be continuous, locally Lipschitz in \((0,\infty)^2\), and zero outside a compact subset of \((0,\infty)\times[0,\infty)\), and assume that \(\int_0^\infty\int_0^\infty\rho^{N-1}y^\alpha|\nabla_{(\rho,y)}h|^2\dd\rho\dd y<\infty\). Then
    \begin{multline}\label{eq:reduced-GSR}
      \int_0^\infty\int_0^\infty
      \rho^{N-1}y^\alpha
      \pt{
        |\nabla_{(\rho,y)}h|^2
        +\frac{(1-\beta_2)\lam}{\rho^2}\abs{h}^2
        -\frac{\beta_1^2}{4(\rho^2+y^2)}\abs{h}^2
      }\dd\rho\dd y
      \\
      =
      \mathcal K(\mu)\int_0^\infty\rho^{N+\alpha-2}\abs{h(\rho,0)}^2\dd\rho
      +
      \int_0^\infty\int_0^\infty\rho^{N-1}y^\alpha w^2\abs{\nabla_{(\rho,y)}\pt{\frac{h}{w}}}^2\dd\rho\dd y.
    \end{multline}
  \end{lemma}

  \begin{proof}[Proof of Lemma~\ref{lem:reduced-GSR}]
    Write \(t=y/\rho\), \(\Phi=(1+t^2)^{\gamma/2}\Vadm\), and
    \(w=\rho^s\Phi(t)\). By \cref{prop:explicit-Kred}, \(w\) is positive,
    smooth for \(y>0\), and continuous up to the boundary, with
    \(w(\rho,0)=\rho^s\). The computations of~\cref{sec:selfsimilar-reduction} show that it solves~\eqref{eq:mode-equation}. Its flux satisfies
    \begin{equation}\label{eq:w-flux}
      \rho^{N-1}y^\alpha\partial_yw(\rho,y)
      =\rho^{s+N+\alpha-2}t^\alpha\Phi'(t)
      \xrightarrow[y\to0^+]{}
      -\mathcal K(\mu)\rho^{N+\alpha-2}w(\rho,0),
    \end{equation}
    uniformly for \(\rho\) in compact subsets of \((0,\infty)\). Indeed, on such a set \(t=y/\rho\to0\) uniformly as \(y\to0^+\), and \(t^\alpha\Phi'(t)=t^\alpha\pt{\Vadm'(t)+\frac{\gamma t}{1+t^2}\Vadm(t)}(1+t^2)^{\gamma/2}\to-K(\Vadm)=-\mathcal K(\mu)\) as \(t\to0^+\) by \cref{prop:explicit-Kred}\ref{it:Vadm-K}.

    Choose \(0<\rho_0<\rho_1\) and \(Y>0\) so that \(h\) is supported
    in \([\rho_0,\rho_1]\times[0,Y]\), and put \(v=h/w\). The product rule gives
    \[
      |\nabla h|^2=w^2|\nabla v|^2+
      \nabla w\cdot\nabla\left(\frac{|h|^2}{w}\right).
    \]
    Set \(\varpi(\rho,y):=\rho^{N-1}y^\alpha\) and \(W(\rho,y):=\beta_1^2/(4(\rho^2+y^2))-(1-\beta_2)\lam/\rho^2\).
    On \((0,\infty)\times(\varepsilon,\infty)\),
    \(|h|^2/w\) is a compactly supported Lipschitz test function.
    Integration by parts, with outward normal \(-e_y\) at the lower boundary, and \(-\diver(\varpi\nabla w)=\varpi Ww\), yield
    \begin{align*}
      & \int_\varepsilon^\infty\!\!\int_0^\infty\varpi\bigl(|\nabla h|^2-W|h|^2\bigr)\dd\rho\dd y
      =\int_\varepsilon^\infty\!\!\int_0^\infty\varpi w^2|\nabla v|^2\dd\rho\dd y\\
      &\hspace{14em}-\int_0^\infty\rho^{N-1}\varepsilon^\alpha
      \partial_yw(\rho,\varepsilon)
      \frac{|h(\rho,\varepsilon)|^2}{w(\rho,\varepsilon)}\dd\rho.
    \end{align*}
    As \(\varepsilon\to0^+\), the gradient integral converges by its
    assumed finiteness. The potential integral converges by dominated convergence, since \(h\) is bounded, \(\rho\) stays away from zero, and \(y^\alpha\) is integrable. Uniform flux convergence and continuity of \(h/w\) give the boundary term in \eqref{eq:reduced-GSR}. The nonnegative remainder converges monotonically; its limit is finite because all other terms have finite limits. This proves the identity.
  \end{proof}

  \begin{remark}\label{rem:GSR-hypotheses}
    The assumptions on \(h\) in \cref{lem:reduced-GSR} hold for every \(h\in C^1_c\pt{(0,\infty)\times[0,\infty)}\), and for every \(h=gw\) with \(g:(0,\infty)\times[0,\infty)\to\C\) Lipschitz and zero outside a compact subset of \((0,\infty)\times[0,\infty)\). Indeed, let \(0<\rho_0<\rho_1\), \(Y>0\) and \(Q_0:=[\rho_0,\rho_1]\times(0,Y]\). By \eqref{eq:Phi-prime} and \cref{prop:explicit-Kred}\ref{it:Vadm-local}, \(t^\alpha\Phi'(t)\) is continuous on \((0,\infty)\) with a finite limit as \(t\to0^+\), so it is bounded for \(t\in(0,Y/\rho_0]\). Hence \(\abs{\partial_yw}=\rho^{s-1}t^{-\alpha}\abs{t^\alpha\Phi'(t)}\leq Cy^{-\alpha}\) on \(Q_0\), and \(\partial_\rho w=\rho^{s-1}\pt{s\Phi-t^{1-\alpha}\,t^\alpha\Phi'}\) is bounded on \(Q_0\), since \(\alpha<1\). Thus \(y^\alpha\abs{\nabla w}^2\leq C\pt{y^\alpha+y^{-\alpha}}\) on \(Q_0\), which is integrable because \(\abs{\alpha}<1\), and \(\abs{\nabla(gw)}^2\leq2\abs{g}^2\abs{\nabla w}^2+2w^2\abs{\nabla g}^2\).
  \end{remark}

  \subsection{The deficit identity and the lower bound}\label{sec:deficit}
  In this subsection we prove \cref{thm:deficit} and \cref{thm:main}\ref{thm:main-lower}. Let \(q_\a[\psi]=\int_{\Sph}|\nabla_\a\psi|^2\dd\omega\) be the
  closed angular form, and recall the deficit from the introduction,
  \[
    \begin{aligned}
      \mathfrak D[u]:={}&\int_{\R^{N+1}_+}y^\alpha\left(
        |\nabla_\Aext u|^2-\frac{\beta_1^2|u|^2}{4(|x|^2+y^2)}
      -\frac{\beta_2\lam|u|^2}{|x|^2}\right)\dd x\dd y\\
      &\hspace{12em}-\mathcal K(\mu)\int_{\R^N}|x|^{\alpha-1}|u(x,0)|^2\dd x.
    \end{aligned}
  \]

  \Cref{prop:full-deficit,prop:energy-norm} below prove \cref{thm:deficit}; there \eqref{eq:deficit-intro} is written in polar coordinates \(x=\rho\omega\).

  \begin{proposition}[Exact deficit identity]\label{prop:full-deficit}
    For every \(u\in\mathcal C_\Aext\),
    \begin{equation}\label{eq:full-deficit}
      \begin{aligned}
        \mathfrak D[u]={}&
        \int_{\Sph}\int_0^\infty\int_0^\infty
        \rho^{N-1}y^\alpha w^2
        \left|\nabla_{(\rho,y)}\left(\frac{u(\rho\omega,y)}{w}\right)\right|^2
        \dd\rho\dd y\dd\omega\\
        &\hspace{4em}+\int_0^\infty\int_0^\infty\rho^{N-3}y^\alpha
        \left(q_\a[u(\rho\,\cdot,y)]-\lam\|u(\rho\,\cdot,y)\|_2^2\right)
        \dd\rho\dd y.
      \end{aligned}
    \end{equation}
    Both terms are nonnegative. In particular, all integrals are finite and
    \begin{equation}\label{eq:step1-inequality}
      \begin{aligned}
        &\mathcal K(\mu)\int_{\R^N}\frac{|u(x,0)|^2}{|x|^{1-\alpha}}\dd x
        +\frac{\beta_1^2}{4}\int_{\R^{N+1}_+}
        \frac{y^\alpha|u|^2}{|x|^2+y^2}\dd x\dd y\\
        &\hspace{10em}+\beta_2\,\lam\,\int_{\R^{N+1}_+}\frac{y^\alpha|u|^2}{|x|^2}\dd x\dd y
        \leq\int_{\R^{N+1}_+}y^\alpha|\nabla_\Aext u|^2\dd x\dd y,
      \end{aligned}
    \end{equation}
    so \(\Kstar\geq\mathcal K(\mu)\).
  \end{proposition}
  \begin{proof}[Proof of Proposition~\ref{prop:full-deficit}]
    The support of \(u\) stays in a compact radial interval away from zero,
    at bounded height, and in a compact angular subset of \(\Sph\setminus\Sigma\).
    Thus all the original integrals are finite. Polar decomposition
    \eqref{eq:grad-spherical} separates the radial--vertical energy from
    \(q_\a\). Add and subtract \(\lam\|u(\rho\,\cdot,y)\|_2^2\) in the latter
    and apply \cref{lem:reduced-GSR} to \(h_\omega(\rho,y)=u(\rho\omega,y)\).
    Integration over \(\omega\) gives \eqref{eq:full-deficit}. The second term
    is nonnegative because
    \begin{equation}\label{eq:angular-lower-bound}
      q_\a[\psi]\geq\lam\|\psi\|_2^2,\quad\text{ for }\;\psi\in H^1_\a(\Sph).
    \end{equation}
  \end{proof}
  To extend \eqref{eq:full-deficit} beyond \(\mathcal C_\Aext\), define the norm and its completion by
  \[
    \begin{aligned}
      \|u\|_*^2&:=\int_{\R^{N+1}_+}y^\alpha
      \left(|\nabla_\Aext u|^2+\frac{|u|^2}{|x|^2}\right)\dd x\dd y
      +\int_{\R^N}|x|^{\alpha-1}|u(x,0)|^2\dd x,\\
      \mathcal D_\Aext&:=\overline{\mathcal C_\Aext}^{\|\cdot\|_*}.
    \end{aligned}
  \]

  \begin{proposition}[Energy-space extension and nonattainment]\label{prop:energy-norm}
    Under \eqref{hyp:lambda}, the norms
    \[
      \|u\|_*\qquad\text{and}\qquad
      \left(\int_{\R^{N+1}_+}y^\alpha|\nabla_\Aext u|^2\dd x\dd y\right)^{1/2}
    \]
    are equivalent on \(\mathcal C_\Aext\). Identity \eqref{eq:full-deficit} extends to \(\mathcal D_\Aext\), and equality in \eqref{eq:step1-inequality} holds there only for \(u=0\). In particular, its optimal trace quotient is not attained in \(\mathcal D_\Aext\).
  \end{proposition}
  \begin{proof}[Proof of Proposition~\ref{prop:energy-norm}]
    Integrating \eqref{eq:hardy-magnetic} in \(y\), and applying
    \cref{prop:full-deficit} with both interior coefficients zero, gives
    \[
      \begin{aligned}
        c_H\int_{\R^{N+1}_+}\frac{y^\alpha|u|^2}{|x|^2}\dd x\dd y
        &\leq\int_{\R^{N+1}_+}y^\alpha|\nabla_\Aext u|^2\dd x\dd y,\\
        K_{N,\a,\alpha}\int_{\R^N}|x|^{\alpha-1}|u(x,0)|^2\dd x
        &\leq\int_{\R^{N+1}_+}y^\alpha|\nabla_\Aext u|^2\dd x\dd y.
      \end{aligned}
    \]
    Here \(c_H=(N-2)^2/4+\lam>0\). By \cref{prop:full-deficit}, \(K_{N,\a,\alpha}\geq\mathcal K(\mu)>0\). Hence
    \[
      \int_{\R^{N+1}_+}y^\alpha|\nabla_\Aext u|^2\dd x\dd y\leq\|u\|_*^2
      \leq(1+c_H^{-1}+K_{N,\a,\alpha}^{-1})\int_{\R^{N+1}_+}y^\alpha|\nabla_\Aext u|^2\dd x\dd y.
    \]
    Next we identify the elements of \(\mathcal D_\Aext\) with functions. Fix a compact set \(F\subset\R^N\setminus\Gamma_\Sigma\) and \(Y>0\). On \(F\), \(\abs x\) is bounded above and below by positive constants and \(x/\abs x\) stays in a compact subset of \(\Sph\setminus\Sigma\), so \(\abs{\A}\leq C_F\) and \(\abs x\leq C_F\) on \(F\). Since \(\nabla_{x,y}u=\nabla_\Aext u+\ii\Aext u\) and \(\Aext=(\A,0)\), every \(u\in\mathcal C_\Aext\) satisfies
    \begin{equation}\label{eq:local-H1}
      \int_0^Y\!\!\int_F y^\alpha\pt{\abs{\nabla_{x,y}u}^2+\abs u^2}\dd x\dd y
      \leq C_F'\|u\|_*^2.
    \end{equation}
    Moreover, for \(v\in\mathcal C_\Aext\),
    \begin{equation}\label{eq:local-trace}
      \|v(\cdot,0)\|_{L^2(F)}^2
      \leq C_{\alpha,Y}\int_0^Y y^\alpha
      \left(\|v(\cdot,y)\|_{L^2(F)}^2+
      \|\partial_yv(\cdot,y)\|_{L^2(F)}^2\right)\dd y.
    \end{equation}
    Indeed, write \(v(0)=v(t)-\int_0^t v'(y)\dd y\), use Cauchy--Schwarz with \(\int_0^t y^{-\alpha}\dd y<\infty\), and average over \(t\in[Y/2,Y]\). Let \((u_k)\subset\mathcal C_\Aext\) be a Cauchy sequence for \(\|\cdot\|_*\). By \eqref{eq:local-H1}, \((u_k)\) converges in \(L^2(F\times(0,Y),y^\alpha\dd x\dd y)\), together with its gradient, to a function \(u\); since \(F\) and \(Y\) are arbitrary and \(\Gamma_\Sigma\times[0,\infty)\) has measure zero, this defines \(u\) almost everywhere in \(\R^{N+1}_+\). By \eqref{eq:local-trace}, \(u_k(\cdot,0)\to u(\cdot,0)\) in \(L^2(F)\), so the limit of \(u_k(\cdot,0)\) in \(L^2(\R^N,\abs x^{\alpha-1}\dd x)\) is \(u(\cdot,0)\). Similarly, the limits of \(y^{\alpha/2}\nabla_\Aext u_k\) and \(y^{\alpha/2}u_k/\abs x\) in \(L^2(\R^{N+1}_+)\) are \(y^{\alpha/2}\nabla_\Aext u\) and \(y^{\alpha/2}u/\abs x\), respectively. In particular, if \(u=0\), then \(\|u_k\|_*\to0\). Hence each element of \(\mathcal D_\Aext\) is a function \(u\) with \(\|u\|_*<\infty\), given by the formula for \(\|\cdot\|_*\), and with trace \(u(\cdot,0)\in L^2(\R^N,\abs x^{\alpha-1}\dd x)\).

    We now extend \eqref{eq:full-deficit} to \(\mathcal D_\Aext\). Let \(\psi_j\) and \(\lambda_j\) be as in \cref{lem:compact-embedding}, and write \(c_j(\psi):=\langle\psi,\psi_j\rangle_{L^2(\Sph)}\). Since \((-\Delta_\a)\) is the operator associated with \(q_\a\) and \(-\Delta_\a\psi_j=\lambda_j\psi_j\), we have \(\int_{\Sph}\nabla_\a\psi\cdot\overline{\nabla_\a\psi_j}\dd\omega=\lambda_jc_j(\psi)\) for \(\psi\in H_\a^1(\Sph)\). Hence the functions \((1+\lambda_j)^{-1/2}\psi_j\) are orthonormal in \(H_\a^1(\Sph)\), and they are complete there, because an element of \(H_\a^1(\Sph)\) orthogonal to all of them is orthogonal to all \(\psi_j\) in \(L^2(\Sph)\). Parseval's identity in \(H_\a^1(\Sph)\) and in \(L^2(\Sph)\) gives
    \[
      H_\a^1(\Sph)=\Bigl\{\psi\in L^2(\Sph):\; \sum_{j=1}^\infty\lambda_j|c_j(\psi)|^2<\infty\Bigr\},
      \qquad
      q_\a[\psi]=\sum_{j=1}^\infty\lambda_j|c_j(\psi)|^2.
    \]
    Let
    \[
      T\psi:=\sum_{j=1}^\infty(\lambda_j-\lam)^{1/2}c_j(\psi)\psi_j,\quad\text{ for }\;\psi\in H_\a^1(\Sph),
    \]
    so that \(q_\a[\psi]-\lam\|\psi\|_2^2=\|T\psi\|_2^2\). The operator \(T\) is closed: if \(\psi_n\to\psi\) and \(T\psi_n\to h\) in \(L^2(\Sph)\), then \(c_j(h)=(\lambda_j-\lam)^{1/2}c_j(\psi)\) for every \(j\), so \(\psi\in H_\a^1(\Sph)\) and \(h=T\psi\). Denote by \(I_1[u]\) and \(I_2[u]\) the two terms on the right of \eqref{eq:full-deficit}; thus
    \[
      I_2[u]=\int_0^\infty\!\!\int_0^\infty\rho^{N-3}y^\alpha\|Tu(\rho\,\cdot,y)\|_2^2\dd\rho\dd y.
    \]
    Let \(u\in\mathcal D_\Aext\), and let \((u_k)\subset\mathcal C_\Aext\) converge to \(u\) in \(\|\cdot\|_*\). Each of the four integrals in \(\mathfrak D\) is the square of the \(L^2\) norm of \(y^{\alpha/2}\nabla_\Aext u\), of \(y^{\alpha/2}u/(\abs x^2+y^2)^{1/2}\), of \(y^{\alpha/2}u/\abs x\), or of \(\abs x^{(\alpha-1)/2}u(\cdot,0)\), so by the identification above \(\mathfrak D[u_k]\to\mathfrak D[u]\). Since \(I_1\) and \(I_2\) are nonnegative, \eqref{eq:full-deficit} applied to \(u_k-u_l\) gives
    \[
      I_1[u_k-u_l]+I_2[u_k-u_l]=\mathfrak D[u_k-u_l]\leq\int_{\R^{N+1}_+}y^\alpha\abs{\nabla_\Aext(u_k-u_l)}^2\dd x\dd y\leq\|u_k-u_l\|_*^2.
    \]
    Hence \(G_k:=\rho^{(N-1)/2}y^{\alpha/2}w\nabla_{(\rho,y)}(u_k/w)\) is a Cauchy sequence in \(L^2(\Sph\times(0,\infty)^2)\), and \(H_k(\rho,y):=\rho^{(N-3)/2}y^{\alpha/2}Tu_k(\rho\,\cdot,y)\) is a Cauchy sequence in \(L^2((0,\infty)^2;L^2(\Sph))\). Let \(G\) and \(H\) be their limits.
    \begin{enumerate}[label=\(\bullet\), itemsep=0em]
      \item On every compact subset of \((\Sph\setminus\Sigma)\times(0,\infty)^2\), the function \(w\) is smooth and bounded below by a positive constant, and \(u_k\to u\) in \(H^1\) by \eqref{eq:local-H1}. Hence \(\nabla_{(\rho,y)}(u_k/w)\to\nabla_{(\rho,y)}(u/w)\) there, and \(G=\rho^{(N-1)/2}y^{\alpha/2}w\nabla_{(\rho,y)}(u/w)\) almost everywhere. Thus \(I_1[u_k]\to I_1[u]\).
      \item Put \(\tilde u_k(\rho,y):=\rho^{(N-3)/2}y^{\alpha/2}u_k(\rho\,\cdot,y)\) and define \(\tilde u\) in the same way, so that \(H_k=T\tilde u_k\). Since \(\rho^{N-3}\dd\rho\dd\omega=\abs x^{-2}\dd x\), we have \(\tilde u_k\to\tilde u\) in \(L^2((0,\infty)^2;L^2(\Sph))\). Passing to a subsequence along which \(\tilde u_k(\rho,y)\) and \(H_k(\rho,y)\) converge in \(L^2(\Sph)\) for almost every \((\rho,y)\), and using that \(T\) is closed, we obtain \(u(\rho\,\cdot,y)\in H_\a^1(\Sph)\) and \(H=T\tilde u\) almost everywhere. Thus \(I_2[u_k]\to I_2[u]\).
    \end{enumerate}
    Letting \(k\to\infty\) in \(\mathfrak D[u_k]=I_1[u_k]+I_2[u_k]\) gives \(\mathfrak D[u]=I_1[u]+I_2[u]\), which is \eqref{eq:full-deficit} for \(u\).

    Finally, let \(u\in\mathcal D_\Aext\) attain equality in \eqref{eq:step1-inequality}, that is, \(\mathfrak D[u]=0\). Then \(I_1[u]=I_2[u]=0\). From \(I_1[u]=0\), for almost every \(\omega\) the function \(u(\rho\omega,y)/w(\rho,y)\) has zero gradient on the connected set \((0,\infty)^2\), so \(u(\rho\omega,y)=w(\rho,y)\psi(\omega)\) for some function \(\psi\) on \(\Sph\); from \(I_2[u]=0\), \(T\psi=0\). Since \(w(\rho,0)=\rho^{-(N+\alpha-1)/2}\),
    \[
      \int_{\R^N}|x|^{\alpha-1}|u(x,0)|^2\dd x=\|\psi\|_2^2\int_0^\infty\frac{1}{\rho}\dd\rho.
    \]
    Finiteness forces \(\psi=0\), proving the claim.
  \end{proof}

  \subsection{The upper bound}\label{sec:upper-bound}
  In this subsection we prove \cref{thm:main}\ref{thm:main-sharp}. We use test functions of the form \(u(x,y)=h(|x|,y)\varphi(x/|x|)\), where \(h\) is a truncation of the profile \(w\) and \(\varphi\in C_c^\infty(\Sph\setminus\Sigma)\).

  These functions are not necessarily smooth up to \(\{y=0\}\), as \(\Vadm(t)=1+\kappa_{\mathrm{adm}}t^{1-\alpha}+O\pt{t^2}\) as \(t\to 0^+\) by \cref{prop:explicit-Kred}\ref{it:Vadm-local}, so they need not belong to \(\mathcal C_\Aext\). However, in the definition \eqref{eq:Kstar-def} of \(\Kstar\) we may replace \(\mathcal C_\Aext\) by the space \(\mathcal D_\Aext\) of \cref{sec:deficit}. Let \(\mathcal O:=\overline{\R^{N+1}_+}\setminus\pt{\Gamma_\Sigma\times[0,\infty)}\). For a continuous function on \(\mathcal O\) that is locally Lipschitz in its interior, the formula for \(\norm{\cdot}_*^2\) defines a possibly infinite quantity. Every integral in \(\Rr\) is continuous with respect to \(\norm{\cdot}_*\), because \(\abs{x}^2+y^2\geq\abs{x}^2\). Therefore
  \begin{equation}\label{eq:Kstar-closure}
    \Kstar
    =
    \inf\{\Rr[u]:\, u\in\mathcal D_\Aext,\ u(\cdot,0)\not\equiv0\}.
  \end{equation}

  The next lemma shows that the test functions described at the beginning of this subsection belong to \(\mathcal D_\Aext\).

  \begin{lemma}[Admissibility of separated test functions]\label{lem:energy-space}
    Let \(h:(0,\infty)\times[0,\infty)\to\R\) be continuous, locally Lipschitz in \((0,\infty)^2\), and zero outside \([\rho_0,\rho_1]\times[0,Y]\) for some \(0<\rho_0<\rho_1\) and \(Y>0\), and assume that \(\int_0^\infty\int_0^\infty y^\alpha|\nabla_{(\rho,y)}h|^2\dd\rho\dd y<\infty\). Then for every \(\varphi\in C_c^\infty(\Sph\setminus\Sigma)\) the function \(u(\rho\omega,y):=h(\rho,y)\varphi(\omega)\) belongs to \(\mathcal D_\Aext\).
  \end{lemma}

  \begin{proof}[Proof of Lemma~\ref{lem:energy-space}]
    Extend \(h\) by zero to \(\rho\leq0\), and let \(h_{\mathrm{ev}}(\rho,y):=h(\rho,\abs{y})\) for \((\rho,y)\in\R^2\). Then \(h_{\mathrm{ev}}\) is continuous, it is locally Lipschitz outside \(\{y=0\}\), its pointwise gradient satisfies \(\int_{\R^2}\abs{y}^\alpha\abs{\nabla h_{\mathrm{ev}}}^2\dd\rho\dd y<\infty\), and, since \(\abs{y}^{-\alpha}\) is locally integrable, \(\nabla h_{\mathrm{ev}}\in L^1(\R^2)\) by the Cauchy--Schwarz inequality. The pointwise gradient is also the distributional gradient. Indeed, for \(f_0\in C_c^\infty(\R^2)\) and \(\ve>0\), integrating by parts on \(\{\abs{y}>\ve\}\), where \(h_{\mathrm{ev}}\) is Lipschitz on compact sets, gives
    \[
      \int_{\{\abs{y}>\ve\}}h_{\mathrm{ev}}\,\partial_yf_0\dd\rho\dd y
      =
      -\int_{\{\abs{y}>\ve\}}\partial_yh_{\mathrm{ev}}f_0\dd\rho\dd y
      -\int_\R h(\rho,\ve)\pt{f_0(\rho,\ve)-f_0(\rho,-\ve)}\dd\rho,
    \]
    and the last term tends to zero as \(\ve\to0^+\) because \(h\) is bounded; the derivative in \(\rho\) is treated in the same way, without boundary terms. Hence \(h_{\mathrm{ev}}\in W^{1,1}(\R^2)\).

    \begingroup\emergencystretch=2em
    \noindent Mollifying \(h_{\mathrm{ev}}\) and restricting to \((0,\infty)\times[0,\infty)\) gives functions \(h_k\in C_c^\infty\pt{(0,\infty)\times[0,\infty)}\), vanishing outside \([\rho_0/2,2\rho_1]\times[0,2Y]\), such that \(h_k\to h\) uniformly and \[
      \int_0^\infty\int_0^\infty y^\alpha|\nabla_{(\rho,y)}(h_k-h)|^2\dd\rho\dd y\to0;
    \]
    the latter holds because \(\abs{y}^\alpha\in A_2(\R^2)\) \cite{Muc1972,Gra2014}. Then \(u_k(\rho\omega,y):=h_k(\rho,y)\varphi(\omega)\) belongs to \(\mathcal C_\Aext\), and \cref{lem:separated}, applied to \(h_k-h\), shows that \(\norm{u_k-u}_*\to0\), because the powers of \(\rho\) are bounded on \([\rho_0/2,2\rho_1]\) and \(y^\alpha\) is integrable near \(0\).\par\endgroup
  \end{proof}

  \Needspace{1em}
  \begin{proposition}[Sharpness of the trace coefficient]\label{prop:Kstar-le-Kred}
    \(\Kstar\leq\mathcal K(\mu)\).
  \end{proposition}

  \begin{proof}[Proof of Proposition~\ref{prop:Kstar-le-Kred}]
    For \(h\) as in \cref{lem:energy-space}, let
    \[
      Q[h]
      :=
      \int_0^\infty\int_0^\infty
      \rho^{N-1}y^\alpha
      \pt{|\nabla_{(\rho,y)}h|^2+\frac{(1-\beta_2)\lam}{\rho^2}h^2
        -\frac{\beta_1^2}{4(\rho^2+y^2)}h^2
      }\dd\rho\dd y,
    \]
    \[
      T_2[h]
      :=
      \int_0^\infty\rho^{N+\alpha-2}h(\rho,0)^2\dd\rho,
      \qquad
      P[h]
      :=
      \int_0^\infty\int_0^\infty\rho^{N-3}y^\alpha h^2\dd\rho\dd y.
    \]
    These quantities are finite by the energy assumption and the compact support separated from \(\rho=0\).

    Assume that \(T_2[h]>0\), let \(\varphi\in C_c^\infty(\Sph\setminus\Sigma)\) with \(\int_{\Sph}\abs{\varphi}^2\dd\omega=1\), and let \(u(\rho\omega,y):=h(\rho,y)\varphi(\omega)\). By \cref{lem:energy-space}, \(u\in\mathcal D_\Aext\), so \(\Kstar\leq\Rr[u]\) by \eqref{eq:Kstar-closure}, and by \cref{lem:separated},
    \begin{equation}\label{eq:R-product}
      \Rr[u]
      =
      \frac{1}{T_2[h]}
      \left(Q[h]+\pt{\int_{\Sph}\abs{\nabla_{\a}\varphi}^2\dd\omega-\lam}P[h]\right).
    \end{equation}

    Let \(\zeta\in C^\infty(\R)\) with \(0\leq\zeta\leq1\), \(\zeta\equiv1\) on \([-1,1]\) and \(\zeta\equiv0\) outside \((-2,2)\), and for \(m\geq2\) let \(\eta_m(\rho):=\zeta\pt{\log\rho/\log m}\). Then \(\eta_m\in C_c^\infty(0,\infty)\), \(\eta_m\equiv1\) on \([1/m,m]\), \(\eta_m\) vanishes outside \((1/m^2,m^2)\), and, substituting \(\varsigma=\log\rho/\log m\),
    \[
      L_m:=\int_0^\infty\eta_m(\rho)^2\frac{\dd\rho}{\rho}\geq2\log m,
      \qquad
      M_m:=\int_0^\infty\eta_m'(\rho)^2\rho\dd\rho=\frac{1}{\log m}\int_\R\zeta'(\varsigma)^2\dd\varsigma.
    \]
    Fix a nonincreasing \(\theta\in C^\infty([0,\infty))\) with \(0\leq\theta\leq1\), \(\theta\equiv1\) on \([0,1]\) and \(\theta\equiv0\) on \([2,\infty)\), let \(\theta_R(y):=\theta(y/R)\) for \(R\geq1\), and let
    \[
      h_{m,R}(\rho,y):=\eta_m(\rho)\,\theta_R(y)\,w(\rho,y).
    \]
    By \cref{rem:GSR-hypotheses}, \cref{lem:reduced-GSR} applies to \(h_{m,R}\), and since \(h_{m,R}/w=\eta_m\theta_R\), \(\abs{\nabla(\eta_m\theta_R)}^2=\theta_R^2\eta_m'^2+\eta_m^2\theta_R'^2\leq\eta_m'^2+\eta_m^2\theta_R'^2\), and \(\abs{\theta_R'}\leq\norm{\theta'}_\infty/R\), \cref{lem:radial-factor} gives
    \[
      Q[h_{m,R}]
      \leq
      \mathcal K(\mu)\,T_2[h_{m,R}]
      +
      M_m\,P_1[\Vadm]
      +
      \frac{\norm{\theta'}_\infty^2}{R^2}\,P_1[\Vadm]\int_0^\infty\eta_m^2\,\rho\dd\rho.
    \]
    Moreover, by \cref{lem:radial-factor} and since \(\theta_R(0)=1\) and \(\theta_R\leq1\),
    \[
      T_2[h_{m,R}]=L_m\,T[\Vadm]=L_m,
      \qquad
      P[h_{m,R}]\leq L_m\,P_1[\Vadm],
    \]
    and these quantities are finite by \cref{prop:explicit-Kred}\ref{it:Vadm-finite}. The function \(h_{m,R}\) satisfies the assumptions of \cref{lem:energy-space}: it is continuous on \((0,\infty)\times[0,\infty)\), locally Lipschitz in \((0,\infty)^2\), and zero outside \([1/m^2,m^2]\times[0,2R]\), and \(\int_0^\infty\int_0^\infty y^\alpha\abs{\nabla h_{m,R}}^2\dd\rho\dd y<\infty\) by \cref{rem:GSR-hypotheses}, since \(\rho^{N-1}\) is bounded below there. Since also \(T_2[h_{m,R}]=L_m>0\), \eqref{eq:R-product} applies to \(h_{m,R}\).

    To control the remaining angular error, the definition
    \eqref{eq:lambda1}--\eqref{eq:HA1} of \(\lam\) gives, for each
    \(\ve>0\), a function \(\varphi_\ve\in C_c^\infty(\Sph\setminus\Sigma)\) such that
    \[
      \int_{\Sph}\abs{\varphi_\ve}^2\dd\omega=1,
      \qquad
      \int_{\Sph}\abs{\nabla_{\a}\varphi_\ve}^2\dd\omega\leq\lam+\ve.
    \]
    Applying \eqref{eq:R-product} to \(u_{m,R,\ve}(\rho\omega,y):=h_{m,R}(\rho,y)\varphi_\ve(\omega)\), and using the estimates above,
    \[
      \begin{aligned}
        \Kstar
        &\leq\Rr[u_{m,R,\ve}]
        \leq \mathcal K(\mu)+\frac{M_m}{L_m}P_1[\Vadm]\\
        &\hspace{12em}+\frac{\norm{\theta'}_\infty^2}{R^2L_m}\,P_1[\Vadm]\int_0^\infty\eta_m^2\,\rho\dd\rho+\ve P_1[\Vadm].
      \end{aligned}
    \]
    We choose one diagonal sequence for all three errors. The change of
    variables above gives
    \(L_m=(\log m)\int_\R\zeta(\varsigma)^2\dd\varsigma\), so
    \[
      \frac{M_m}{L_m}
      =\frac{\int_\R\zeta'(\varsigma)^2\dd\varsigma}
      {(\log m)^2\int_\R\zeta(\varsigma)^2\dd\varsigma}
      \leq\frac{1}{2(\log m)^2}
      \int_\R\zeta'(\varsigma)^2\dd\varsigma
      \longrightarrow0.
    \]
    Take \(m_n=e^n\), \(R_n=m_n^3\), and \(\ve_n=1/n\). Since \(\int_0^\infty\eta_m^2\,\rho\dd\rho \leq m^4/2\), the three errors are respectively \(O(1/n^2)\), \(O(e^{-2n}/n)\), and \(O(1/n)\). For each \(n\), \cref{lem:energy-space} and the continuity of the quotient in \eqref{eq:Kstar-closure} give \(\tilde u_n\in\mathcal C_\Aext\) whose quotient differs from that of \(u_{m_n,R_n,\ve_n}\) by less than \(1/n\). Hence \(\Rr[\tilde u_n]\to\mathcal K(\mu)\).
  \end{proof}

  \section{Comparison with the nonmagnetic case}\label{sec:non-magnetic}
  In this section we compare \(\mathcal K(\mu)\) with the constant of the nonmagnetic problem, and we characterize, in terms of capacity and gauges, when \(\lam=0\).

  When \(\lam=0\), which by \eqref{hyp:lambda} forces \(N\geq3\) and is characterized in \cref{prop:lambda1-zero} below, we have \(\mu=N-2\), and
  \begin{equation}\label{eq:K-explicit}
    \mathcal K(N-2)
    =(1-\alpha)
    \frac{
      \Gamma\pt{\frac{1+\alpha}{2}}
      \Gamma\pt{\frac{N+1-\alpha-\beta_1}{4}}
      \Gamma\pt{\frac{N+1-\alpha+\beta_1}{4}}
    }{
      \Gamma\pt{\frac{3-\alpha}{2}}
      \Gamma\pt{\frac{N+\alpha-1-\beta_1}{4}}
      \Gamma\pt{\frac{N+\alpha-1+\beta_1}{4}}
    },
  \end{equation}
  which is the sharp nonmagnetic constant \(H_{N,\alpha,\beta_1}\)
  in \eqref{eq:full-nonmag-inequality}; in the notation of
  \cite[Theorem 1]{Tzi2015}, set \(n=N+1\) and \(b=\beta_1+2-\alpha\). The condition \(2-\alpha\leq b<n\) is then exactly
  \(0\leq\beta_1<N+\alpha-1\). The Euclidean specialization of
  \cite[Theorem 1]{Tzi2022} gives the same coefficient, with \(n=N\)
  and the same parameter \(b=\beta_1+2-\alpha\).
  For \(\beta_1=0\), using \(\Gamma((3-\alpha)/2)=((1-\alpha)/2)\Gamma((1-\alpha)/2)\), the constant becomes
  \begin{equation}\label{eq:nonmagnetic-H}
    \mathcal K(N-2,\alpha,0)
    = 2 \frac{\Gamma\pt{\frac{1+\alpha}{2}}}{\Gamma\pt{\frac{1-\alpha}{2}}}\frac{\Gamma\pt{\frac{N+1-\alpha}{4}}^2}{\Gamma\pt{\frac{N+\alpha-1}{4}}^2}
    = H_{N,\alpha},
  \end{equation}
  the sharp constant of Tzirakis in \eqref{eq:trace-hardy}.
  \par\Needspace{9\baselineskip}
  \begin{proposition}[Strict magnetic improvement]\label{prop:comparison}
    Let \(\beta_1=\beta_2=0\), and assume \(\lam>0\). Then
    \[
      K_{N,\a,\alpha}=\mathcal K(\mu)>\mathcal K(N-2)=H_{N,\alpha}.
    \]
  \end{proposition}

  \begin{proof}[Proof of Proposition~\ref{prop:comparison}]
    By \cref{thm:main}\ref{thm:main-sharp}, \(K_{N,\a,\alpha}=\mathcal K(\mu)\) with \(\mu=\sqrt{(N-2)^2+4\lam}\), and \(\mu>N-2\geq0\) because \(\lam>0\). The strict inequality follows from \cref{lem:K-increasing} with \(\beta_1=0\), and the last equality is \eqref{eq:nonmagnetic-H}.
  \end{proof}
  For comparison on the same angular domain, let
  \(\lambda_{1,0}^{\Sigma}\) be the first eigenvalue of the closure
  of \(\psi\mapsto\int_{\Sph}|\nabla_{\Sph}\psi|^2\dd\omega\)
  on \(C_c^\infty(\Sph\setminus\Sigma)\). The diamagnetic
  inequality and approximation in the form norms give
  \(\lam\geq\lambda_{1,0}^{\Sigma}\). If \(\Sigma\) has
  positive \(H^1\)-capacity, this closure may impose a Dirichlet
  condition on \(\Sigma\), as in the following two examples. Both rely on the ground state identity: if \(G\subset\Sph\) is open, \(f\in C^\infty(G)\) is positive, and \(-\Delta_{\Sph}f=\lambda f\) in \(G\), then every \(\psi\in C_c^\infty(G)\) satisfies, with \(v=\psi/f\),
  \begin{equation}\label{eq:GSR-sphere}
    \int_G\abs{\nabla_{\Sph}\psi}^2\dd\omega
    =\int_G f^2\abs{\nabla_{\Sph}v}^2\dd\omega+\lambda\int_G\abs{\psi}^2\dd\omega.
  \end{equation}
  Indeed, \(\abs{\nabla(fv)}^2=f^2\abs{\nabla v}^2+\nabla f\cdot\nabla(f\abs v^2)\), and integrating the last term by parts gives \(\int_G(-\Delta_{\Sph}f)f\abs v^2\dd\omega=\lambda\int_G\abs\psi^2\dd\omega\).
  \begin{enumerate}[label=\(\bullet\), itemsep=0em]
    \item Let \(N=3\), \(\a=0\), and let \(\Sigma=\{\omega_3=0\}\) be the equator, so that \(\Ss^2\setminus\Sigma\) is the union of the open hemispheres \(G_\pm=\{\pm\omega_3>0\}\). Every function in \(C_c^\infty(\Ss^2\setminus\Sigma)\) is the sum of its restrictions to \(G_+\) and \(G_-\), so it suffices to consider \(\psi\in C_c^\infty(G_\pm)\). The function \(f=\pm\omega_3\) is positive in \(G_\pm\), and it is the restriction of a spherical harmonic of degree one, so \(-\Delta_{\Ss^2}f=2f\). By \eqref{eq:GSR-sphere}, \(\int\abs{\nabla_{\Ss^2}\psi}^2\geq2\int\abs\psi^2\), hence \(\lambda_{1,0}^{\Sigma}\geq2\). Conversely, \(f\chi_{G_\pm}\) is Lipschitz on \(\Ss^2\) and vanishes on \(\Sigma\), so it belongs to the closure, and \(\int_{G_\pm}\abs{\nabla f}^2=\int_{G_\pm}f(-\Delta f)=2\int_{G_\pm}f^2\) because \(f=0\) on \(\partial G_\pm\). Thus \(\lambda_{1,0}^{\Sigma}=2\), whereas the first eigenvalue of \(-\Delta_{\Ss^2}\) is \(0\).
    \item Let \(N=2\) and \(\Sigma=\{(1,0)\}\), so that \(\Ss^1\setminus\Sigma\) is the interval \(0<\theta<2\pi\). The function \(f=\sin(\theta/2)\) is positive there and satisfies \(-f''=f/4\). By \eqref{eq:GSR-sphere}, \(\lambda_{1,0}^{\Sigma}\geq1/4\), and \(f\), which belongs to \(H_0^1(0,2\pi)\), gives equality. Thus \(\lambda_{1,0}^{\Sigma}=1/4\), whereas the first eigenvalue of \(-\partial_\theta^2\) on \(\Ss^1\) is \(0\). Moreover, \(\lam=\lambda_{1,0}^{\Sigma}=1/4\) for every \(\a\) (\cref{sec:example-2D}).
  \end{enumerate}
  Thus the increase over \(H_{N,\alpha}\) can be purely geometric, and
  strict magnetic improvement relative to the same domain
  requires \(\lam>\lambda_{1,0}^{\Sigma}\).

  We now characterize when \(\lam=0\). All gradients below are angular, and \(\Omega=\Sph\setminus\Sigma\). We identify \(\a\) with its metric-dual one-form, and \(d\a\) denotes its exterior derivative. A gauge is a smooth map
  \(g:\Omega\to\C\) with \(|g|=1\); it need not have a globally
  defined real phase. When \(\nabla g=\ii\a\, g\), multiplication by
  \(g\) satisfies \(\nabla_\a(gv)=g\nabla v\), so it maps \(C_c^\infty(\Omega)\) onto itself and extends to an isometry from the closure of \(C_c^\infty(\Omega)\) in \(H^1(\Sph)\) onto \(H_\a^1(\Sph)\).

  For \(\Sigma=\emptyset\) and \(N\geq3\), the regular equality
  criterion is contained in \cite[Lemma A.2]{FeFeTe2011}.
  \begin{proposition}[Vanishing of the lowest angular eigenvalue]\label{prop:lambda1-zero}
    Under the angular assumptions of \cref{thm:main}, for any \(N\geq2\),
    \[
      \lam=0
      \quad\Longleftrightarrow\quad
      \operatorname{Cap}_{H^1}(\Sigma)=0
      \ \text{and}\
      \exists\,g\in C^\infty(\Omega;\C),\quad |g|=1,\quad\nabla g=\ii\a\,g.
    \]
    In particular, \(\lambda_{1,0}^{\Sigma}=0\) exactly when \(\Sigma\) has
    zero capacity. When \(\Sigma=\emptyset\) and \(N\geq3\), the existence of a smooth \(g\) with \(\abs g=1\) and \(\nabla g=\ii\a\,g\) is equivalent to \(\a=\nabla\chi\) for a smooth real function
    \(\chi\), and also to \(d\a=0\); hence \(\lam=0\) if and only if \(d\a=0\).
  \end{proposition}
  \begin{proof}[Proof of Proposition~\ref{prop:lambda1-zero}]
    Suppose \(\lam=0\), and let \(\psi\) be a normalized first eigenfunction.
    Choose \((\psi_n)\subset C_c^\infty(\Omega)\) with \(\psi_n\to\psi\) in \(H_\a^1(\Sph)\). Diamagnetism gives \(|\psi|\in H^1(\Sph)\) with zero weak gradient; this follows first by weak compactness of \(|\psi_n|\) and then by the lower semicontinuity of the nonmagnetic energy. Since \(\Sph\) is connected, \(|\psi|=c_*>0\) almost everywhere. Local elliptic regularity and \(\nabla_\a\psi=0\) give a smooth representative on \(\Omega\), so \(g=\psi/c_*\) is a parallel unit gauge there. The globally Lipschitz functions
    \[
      \chi_n=1-\min\{1,|\psi_n|/c_*\}
    \]
    equal one near \(\Sigma\). Their \(L^2\)-norms tend to zero because
    \(|\psi_n|\to c_*\) in \(L^2\), and \(\|\nabla\chi_n\|_2\leq c_*^{-1}\|\nabla_\a\psi_n\|_2\to0\). They prove that \(\Sigma\) has zero capacity.

    Conversely, if the two stated conditions hold, then \(g\) is bounded and has zero magnetic energy. By~\cref{lem:capacity-bounded}, \(g\in H^1_\a(\Sph)\), and hence \(\lam=0\). Taking \(\a=0\) gives the nonmagnetic assertion. On a simply connected sphere \(\Ss^{N-1}\), \(N\geq3\), a smooth circle-valued gauge has a smooth real lift; closed one-forms are exact. This proves the last equivalences.
  \end{proof}

  \section{The fractional magnetic Hardy inequality}\label{sec:fractional}

  Throughout, \(\sigma\in(0,1)\), \(\alpha=1-2\sigma\), \(\beta_1=\beta_2=0\), and we assume~\eqref{hyp:lambda}, so that \(\mu=\sqrt{(N-2)^2+4\lam}>0\). Let \(L:=-\Delta_\A\) be the nonnegative self-adjoint operator defined in the introduction. Every \(f\in C_c^\infty(\R^N\setminus\Gamma_\Sigma)\) belongs to the domain \(D(L)\) of \(L\), so \(\int_{\R^N}f\,\overline{L^\sigma f}\dd x=\|L^{\sigma/2}f\|_{L^2}^2\).

  The form domain of \(L\) is the closure \(V\) of \(C_c^\infty(\R^N\setminus\Gamma_\Sigma)\) with respect to the norm
  \[
    \norm{h}_V:=\pt{\int_{\R^N}\pt{\abs{\nabla_\A h}^2+\abs{h}^2}\dd x}^{1/2},
  \]
  and the closed form is \(q(h):=\int_{\R^N}|\nabla_\A h|^2\dd x\), for \(h\in V\). Here \(\nabla_\A h=\nabla h-\ii \A h\) is understood distributionally on \(\R^N\setminus\Gamma_\Sigma\), and as the
  \(L^2\)-limit of magnetic gradients of approximants in \(C_c^\infty(\R^N\setminus\Gamma_\Sigma)\). Since \(\Gamma_\Sigma\) has measure zero, \(C_c^\infty(\R^N\setminus\Gamma_\Sigma)\) is dense in \(L^2(\R^N)\), and, as in the proof of \cref{lem:compact-embedding}, \(q\) is closable on this space because it is the form of the symmetric nonnegative operator \((-\Delta_\A)\).
  We write \(q(h,k)\) for the associated sesquilinear form, so that
  \(q(h,k)=\langle Lh,k\rangle_{L^2}\) for \(h\in D(L)\) and \(k\in V\). Moreover, \(V=D(L^{1/2})\) and \(q(h)=\|L^{1/2}h\|_2^2\) for \(h\in V\).

  Recall that by the spectral theorem \cite[Theorem 13.30]{Rud1973}, there is a unique projection-valued measure \(E\) on the Borel subsets of \([0,\infty)\) such that \(L=\int_{[0,\infty)}\lambda\dd E(\lambda)\). The spectral measure of \(h\in L^2(\R^N)\) with respect to \(L\) is the finite Borel measure
  \[
    \mu_h(I):=\langle E(I)h,h\rangle_{L^2}=\norm{E(I)h}_2^2,\qquad I\subset[0,\infty)\text{ Borel},
  \]
  so that \(\mu_h([0,\infty))=\norm{h}_2^2\). For a Borel function \(\phi\colon[0,\infty)\to\C\), the operator \(\phi(L):=\int_{[0,\infty)}\phi(\lambda)\dd E(\lambda)\) has domain \(D(\phi(L))=\{h\in L^2(\R^N):\int_{[0,\infty)}|\phi|^2\dd\mu_h<\infty\}\), and \cite[Theorem 13.24]{Rud1973}
  \[
    \norm{\phi(L)h}_2^2=\int_{[0,\infty)}|\phi(\lambda)|^2\dd\mu_h(\lambda),\qquad h\in D(\phi(L)).
  \]
  For a Borel set \(I\subset[0,\infty)\), the spectral projection of \(L\) associated with \(I\) is \(\mathbf1_I(L)=E(I)\). Since \(E(J)E(I)=E(J\cap I)\), the spectral measure of \(\mathbf1_I(L)h\) is the restriction of \(\mu_h\) to \(I\). In particular, \(\norm{h}_V^2=\int_{[0,\infty)}(1+\lambda)\dd\mu_h(\lambda)\) for \(h\in V\).

  \subsection{The variational extension principle}
  In this subsection we show that the Stinga--Torrea extension of \(f\) minimizes the weighted energy among all extensions of \(f\), and we compute the minimum.

  We identify a function \(v\) on \(\R^{N+1}_+\) with the map \(y\mapsto v(y):=v(\cdot,y)\in L^2(\R^N)\), and, for maps \(v:(0,\infty)\to L^2(\R^N)\) that are locally absolutely continuous on \((0,\infty)\) and satisfy \(v(y)\in V\) for almost every \(y\), we use the energy integral
  \[
    \mathcal J[v]:=\int_0^\infty y^\alpha\bigl(q(v(y))+\|v'(y)\|_{L^2}^2\bigr)\dd y
    \in[0,\infty].
  \]
  For finite-energy maps, \(v(0)\) denotes the canonical \(L^2\) trace \(\lim_{y\to0^+}v(y)\). This limit exists because, for \(0<z<y\),
  \[
    \|v(y)-v(z)\|_{L^2}
    \leq
    \left(\int_z^y t^\alpha\|v'(t)\|_{L^2}^2\dd t\right)^{1/2}
    \left(\frac{y^{1-\alpha}-z^{1-\alpha}}{1-\alpha}\right)^{1/2},
  \]
  and \(\alpha<1\). We always use this continuous representative at
  \(y=0\). For two finite-energy maps \(v_1\) and \(v_2\), the associated sesquilinear pairing is
  \[
    \mathcal J[v_1,v_2]:=\int_0^\infty y^\alpha\bigl(q(v_1(y),v_2(y))+\langle v_1'(y),v_2'(y)\rangle_{L^2}\bigr)\dd y,
  \]
  so that \(\mathcal J[v]=\mathcal J[v,v]\), and \(\mathcal J[\cdot]^{1/2}\) is a seminorm. For
  \(u\in\mathcal C_\Aext\), \(\mathcal J[u]=\int_{\R^{N+1}_+}y^\alpha|\nabla_\Aext u|^2\dd x\dd y\).
  Since \(\abs{x}^{1-\alpha}=\abs{x}^{2\sigma}\), \cref{thm:main}\ref{thm:main-lower} with \(\beta_1=\beta_2=0\) reads
  \begin{equation}\label{eq:trace-C}
    \mathcal K(\mu)\int_{\R^N}\frac{\abs{u(x,0)}^2}{\abs{x}^{2\sigma}}\dd x
    \leq
    \mathcal J[u],
    \qquad u\in\mathcal C_\Aext.
  \end{equation}
  For \(f\in L^2(\R^N)\), let
  \begin{equation}\label{eq:U-def}
    U_f(\cdot,y)
    :=
    \frac{y^{2\sigma}}{4^\sigma\Gamma(\sigma)}
    \int_0^\infty e^{-tL}f\,e^{-\frac{y^2}{4t}}\,\frac{\dd t}{t^{1+\sigma}},
    \qquad y>0,
  \end{equation}
  and \(U_f(\cdot,0):=f\). By \cite{StTo2010}, if \(f\in C_c^\infty(\R^N\setminus\Gamma_\Sigma)\), then \(U_f\) solves the extension problem
  \begin{equation}\label{eq:ST-extension}
    \partial_y\pt{y^\alpha\partial_yU_f}=y^\alpha LU_f
    \quad\text{in }\R^{N+1}_+,
    \quad
    U_f(\cdot,0)=f,
    \quad
    -\lim_{y\to 0^+}y^\alpha\partial_yU_f(\cdot,y)=\kappa_\sigma L^\sigma f.
  \end{equation}

  \begin{proposition}[Variational extension principle]\label{prop:variational-extension}
    For \(f\in L^2(\R^N)\), let \(\mathcal E_f\) be the set of finite-energy extensions of \(f\) to \(\R^{N+1}_+\), that is, the set of maps \(v:(0,\infty)\to L^2(\R^N)\) such that
    \begin{enumerate}[label=\(\bullet\),itemsep=0em]
      \item \(v\) is locally absolutely continuous on \((0,\infty)\) as a map into \(L^2(\R^N)\), and \(v(y)\in V\) for almost every \(y>0\);
      \item \(v\) has finite energy,
        \(
          \mathcal J[v]<\infty
        \);
      \item the trace of \(v\) defined above is \(v(0)=f\).
    \end{enumerate}
    Then
    \begin{equation}\label{eq:bochner-minimum}
      \begin{gathered}
        \mathcal E_f\ne\emptyset
        \quad\Longleftrightarrow\quad f\in D(L^{\sigma/2}),\\
        \min_{v\in\mathcal E_f}\mathcal J[v]=\kappa_\sigma\|L^{\sigma/2}f\|_2^2.
      \end{gathered}
    \end{equation}
    For such \(f\), the map \(U_f\) of \eqref{eq:U-def} is the unique
    minimizer, and every \(v\in\mathcal E_f\) satisfies
    \begin{equation}\label{eq:bochner-orthogonality}
      \mathcal J[v]-\kappa_\sigma\|L^{\sigma/2}f\|_2^2=\mathcal J[v-U_f],
      \qquad
      \mathcal J[U_f,v-U_f]=0.
    \end{equation}
  \end{proposition}

  \begin{proof}[Proof of Proposition~\ref{prop:variational-extension}]
    The Poisson multiplier is the function
    \begin{equation}\label{eq:poisson-multiplier}
      \Upsilon(s):=\frac1{\Gamma(\sigma)}
      \int_0^\infty e^{-r-s^2/(4r)}r^{\sigma-1}\dd r
      =\frac{2^{1-\sigma}}{\Gamma(\sigma)}\,s^\sigma K_\sigma(s),
      \qquad U_f(y)=\Upsilon(yL^{1/2})f,
    \end{equation}
    where \(K_\sigma\) is the modified Bessel function of the second kind. The identity \(U_f(y)=\Upsilon(yL^{1/2})f\) is the Poisson formula \cite[(1.9)]{StTo2010}, and for each \(\lambda>0\) the function \(y\mapsto\Upsilon(y\lambda^{1/2})\) is the solution of the extension problem computed in \cite[Section 3.1, (3.8)]{StTo2010}. As shown there,
    \[
      (s^{1-2\sigma}\Upsilon')'=s^{1-2\sigma}\Upsilon,
      \qquad
      \Upsilon(0)=1,
      \qquad
      -s^{1-2\sigma}\Upsilon'(s)\xrightarrow[s\to0^+]{}\kappa_\sigma,
    \]
    and \(\Upsilon(s)\) and \(s^{1-2\sigma}\Upsilon'(s)\) are \(O(s^{1/2}e^{-s})\) as \(s\to\infty\). Moreover, by \cite[Vol.~II, \S7.11, (21), and \S7.2.2, (14)]{ErMagOb1953},
    \[
      s^{1-2\sigma}\Upsilon'(s)=-\frac{2^{1-\sigma}}{\Gamma(\sigma)}\,s^{1-\sigma}K_{1-\sigma}(s),
    \]
    and \(K_\nu>0\) on \((0,\infty)\) by \cite[Vol.~II, \S7.12, (21)]{ErMagOb1953}, so \(\Upsilon'<0\) and \(0<\Upsilon\leq1\). Near \(s=0\), \(\Upsilon\) is bounded and \(s^{1-2\sigma}\Upsilon'(s)^2=s^{2\sigma-1}\pt{s^{1-2\sigma}\Upsilon'(s)}^2=O(s^{2\sigma-1})\), which is integrable. Hence, multiplying \((s^{1-2\sigma}\Upsilon')'=s^{1-2\sigma}\Upsilon\) by \(\Upsilon\) and integrating by parts on \((0,\infty)\), the boundary term at infinity vanishes, and \(\Upsilon(0)=1\) together with \(-s^{1-2\sigma}\Upsilon'(s)\to\kappa_\sigma\) gives
    \[
      \int_0^\infty s^{1-2\sigma}\pt{\Upsilon(s)^2+\Upsilon'(s)^2}\dd s
      =-\lim_{s\to0^+}s^{1-2\sigma}\Upsilon'(s)\Upsilon(s)
      =\kappa_\sigma.
    \]
    The change of variables \(s=y\lambda^{1/2}\), with \(\alpha=1-2\sigma\), gives
    \begin{equation}\label{eq:scalar-poisson-energy}
      \int_0^\infty y^\alpha
      \left(\lambda\Upsilon(y\lambda^{1/2})^2+
      \left|\partial_y\Upsilon(y\lambda^{1/2})\right|^2\right)\dd y
      =\kappa_\sigma\lambda^\sigma,\qquad \lambda>0.
    \end{equation}

    To justify the comparison for an arbitrary finite-energy competitor,
    we first restrict to bounded spectral intervals away from zero. By \eqref{hyp:lambda}
    and the closed-form version of \eqref{eq:hardy-magnetic},
    \(\ker L=\{0\}\). Let \(P_m=\mathbf1_{[1/m,m]}(L)\), \(m\geq2\).
    For \(v\in\mathcal E_f\), put \(v_m=P_mv\), \(f_m=P_mf\), and
    \(U_m:=U_{f_m}=\Upsilon(yL^{1/2})f_m\). On this spectral band the extension
    equation holds in \(L^2\), and its flux at zero is
    \(-\kappa_\sigma L^\sigma f_m\). Formula~\eqref{eq:scalar-poisson-energy} and the spectral theorem show that
    \[
      \mathcal J[U_m]=\kappa_\sigma\|L^{\sigma/2}f_m\|_2^2.
    \]

    Since \((s^{1-2\sigma}\Upsilon')'=s^{1-2\sigma}\Upsilon\), the band-limited map \(U_m\) satisfies \((y^\alpha U_m')'=y^\alpha LU_m\) in \(L^2(\R^N)\), with \(U_m(y)\in D(L)\). As \(v_m-U_m\) is locally absolutely continuous, integration by parts against \(v_m-U_m\) on \((\delta,Y)\), \(0<\delta<Y\), gives
    \[
      \begin{aligned}
        &\int_{\delta}^{Y}y^\alpha\bigl(q(U_m(y),v_m(y)-U_m(y))+\langle U_m'(y),v_m'(y)-U_m'(y)\rangle_2\bigr)\dd y\\
        &\quad=\bigl\langle y^\alpha U_m'(y),v_m(y)-U_m(y)\bigr\rangle_2\Big|_{\delta}^{Y}.
      \end{aligned}
    \]
    We show that both boundary terms tend to zero. By the properties of \(\Upsilon\) listed above, \(M:=\sup_{s>0}\abs{s^{1-2\sigma}\Upsilon'(s)}<\infty\), and
    \[
      y^\alpha U_m'(y)=L^\sigma\phi_y(L)f_m,
      \qquad
      \phi_y(\lambda):=(y\lambda^{1/2})^{1-2\sigma}\Upsilon'(y\lambda^{1/2}).
    \]
    \begin{enumerate}[label=\(\bullet\),itemsep=0em]
      \item As \(\delta\to0^+\): \(\|\delta^\alpha U_m'(\delta)\|_2\leq M\|L^\sigma f_m\|_2\), while \(v_m(\delta)-U_m(\delta)\to f_m-f_m=0\) in \(L^2(\R^N)\), because the trace of \(v_m=P_mv\) is \(P_mf=f_m\) and \(\Upsilon(0)=1\). Hence the term at \(\delta\) tends to zero.
      \item As \(Y\to\infty\): the spectral measure of \(f_m\) is supported in \([1/m,m]\), where \(Y\lambda^{1/2}\geq Y/\sqrt m\). Since the function \(s^{1-2\sigma}\Upsilon'(s)\) decays like \(s^{1/2}e^{-s}\), this gives \(\|Y^\alpha U_m'(Y)\|_2\leq C_m(1+Y)^{1/2}e^{-Y/\sqrt m}\). On the other hand, \(\|U_m(Y)\|_2\leq\|f_m\|_2\) because \(0<\Upsilon\leq1\), and for \(Y\geq1\), by Cauchy--Schwarz and \((1-\alpha)/2=\sigma\),
        \[
          \begin{aligned}
            \|v_m(Y)-v_m(1)\|_2
            &\leq\pt{\int_1^Y t^{-\alpha}\dd t}^{1/2}\pt{\int_1^Y t^\alpha\|v_m'(t)\|_2^2\dd t}^{1/2}\\
            &\leq C\,Y^{\sigma}\pt{\int_0^\infty t^\alpha\|v'(t)\|_2^2\dd t}^{1/2}.
          \end{aligned}
        \]
        Hence the term at \(Y\) is at most \(C_m(1+Y)^{1/2+\sigma}e^{-Y/\sqrt m}\), which tends to zero.
    \end{enumerate}
    The energy pairing is integrable by Cauchy--Schwarz, so
    \[
      \mathcal J[v]\geq\mathcal J[v_m]=\kappa_\sigma\|L^{\sigma/2}f_m\|_2^2+\mathcal J[v_m-U_m].
    \]
    Since the spectral measure of \(f_m=P_mf\) is the restriction of \(\mu_f\) to \([1/m,m]\),
    \[
      \|L^{\sigma/2}f_m\|_2^2=\int_{[1/m,m]}\lambda^\sigma\dd\mu_f(\lambda).
    \]
    The intervals \([1/m,m]\) increase to \((0,\infty)\), and \(\lambda^\sigma=0\) at \(\lambda=0\), so the monotone convergence theorem gives
    \[
      \kappa_\sigma\int_{[0,\infty)}\lambda^\sigma\dd\mu_f(\lambda)
      =\lim_{m\to\infty}\kappa_\sigma\|L^{\sigma/2}f_m\|_2^2
      \leq\mathcal J[v]<\infty.
    \]
    Since \(D(L^{\sigma/2})=\{f\in L^2(\R^N):\int_{[0,\infty)}\lambda^\sigma\dd\mu_f<\infty\}\) and \(\|L^{\sigma/2}f\|_2^2=\int_{[0,\infty)}\lambda^\sigma\dd\mu_f\)~\cite[Section 2]{StTo2010}, this proves \(f\in D(L^{\sigma/2})\) and the lower bound in~\eqref{eq:bochner-minimum}.

    Conversely, let \(f\in D(L^{\sigma/2})\). Since \(0<\Upsilon\leq1\) and \(\Upsilon(y\lambda^{1/2})\to1\) as \(y\to0^+\), dominated convergence for the spectral measure \(\mu_f\) shows that \(U_f(y)\to f\) in \(L^2(\R^N)\), so \(U_f\) has trace \(f\). By the behavior of \(\Upsilon\) and \(s^{1-2\sigma}\Upsilon'\) at zero and at infinity listed above, \(s\Upsilon(s)\) and \(s\Upsilon'(s)=s^{2\sigma}\cdot s^{1-2\sigma}\Upsilon'(s)\) are bounded on \([0,\infty)\). With \(s=y\lambda^{1/2}\), this gives, for \(y>0\),
    \[
      \sup_{\lambda\geq0}\lambda\,\Upsilon(y\lambda^{1/2})^2\leq\frac{C}{y^2},
      \qquad
      \sup_{\lambda\geq0}\abs{\partial_y\Upsilon(y\lambda^{1/2})}=\sup_{\lambda\geq0}\abs{\lambda^{1/2}\Upsilon'(y\lambda^{1/2})}\leq\frac{C}{y}.
    \]
    By the first bound, \(U_f(y)\in D(L^{1/2})=V\) for every \(y>0\). By the second bound and dominated convergence, \(y\mapsto U_f(y)\) is \(C^1\) from \((0,\infty)\) into \(L^2(\R^N)\), with \(U_f'(y)=L^{1/2}\Upsilon'(yL^{1/2})f\); in particular, it is locally absolutely continuous. Tonelli's theorem and~\eqref{eq:scalar-poisson-energy} give
    \[
      \mathcal J[U_f]=\kappa_\sigma\|L^{\sigma/2}f\|_2^2,
    \]
    so \(U_f\in\mathcal E_f\). The change of variables in~\eqref{eq:poisson-multiplier} recovers \eqref{eq:U-def} as a Bochner integral.

    Let \(v\in\mathcal E_f\). Since \(\ker L=\{0\}\), \(P_m=\mathbf1_{[1/m,m]}(L)\) converges strongly to the identity, so \((I-P_m)w\to0\) in \(L^2(\R^N)\) for every \(w\in L^2(\R^N)\), with \(\|(I-P_m)w\|_2\leq\|w\|_2\). As \(P_m\) commutes with \(L^{1/2}\) and with \(\partial_y\),
    \[
      \mathcal J[(I-P_m)v]=\int_0^\infty y^\alpha\bigl(\|(I-P_m)L^{1/2}v(y)\|_2^2+\|(I-P_m)v'(y)\|_2^2\bigr)\dd y,
    \]
    and dominated convergence in \(y\), with dominating function \(y^\alpha(q(v(y))+\|v'(y)\|_2^2)\), shows that this tends to zero. Thus \(\mathcal J[v-P_mv]\to0\), and in the same way \(\mathcal J[U_f-U_m]\to0\), because \(P_mU_f=U_m\). Since \(v_m-U_m=P_m(v-U_f)\), the orthogonality relation obtained above for \(v_m\) and \(U_m\) passes to the limit by Cauchy--Schwarz and gives the second identity in \eqref{eq:bochner-orthogonality}. Expanding the energy of \(v=U_f+(v-U_f)\) and using this orthogonality together with the energy identity for \(U_f\) gives the first one. If the energy of \(v\) equals \(\kappa_\sigma\|L^{\sigma/2}f\|_2^2\), the first identity forces \((v-U_f)'=0\) almost everywhere; since \(v-U_f\) is locally absolutely continuous with zero trace, \(v=U_f\).
  \end{proof}

  The sharp constant in \cref{thm:main} is defined on \(\mathcal C_\Aext\), whereas \cref{prop:variational-extension} concerns the class \(\mathcal E_f\). The next two lemmas approximate the extension \(U_f\) in energy by functions in \(\mathcal C_\Aext\).

  \begin{lemma}[Approximation by smooth core functions]\label{lem:approximation}
    Let \(Y>0\), and let \(v:[0,\infty)\to V\) be continuous for the norm of \(V\), with \(v(y)=0\) for \(y\geq Y\). Assume that \(v\) is absolutely continuous on \([0,Y]\) as a map into \(L^2(\R^N)\), with \(\int_0^Yy^\alpha\norm{v'(y)}_{L^2}^2\dd y<\infty\). Then there are \(u_n\in\mathcal C_\Aext\) such that
    \[
      \mathcal J[u_n-v]\longrightarrow0
    \]
    and \(u_n(\cdot,0)\to v(0)\) in \(V\).
  \end{lemma}

  \begin{proof}[Proof of Lemma~\ref{lem:approximation}]
    Extend \(v\) evenly to \(v_{\mathrm{ev}}(y)=v(|y|)\) on \(\mathbb R\).
    This is continuous with compact support as a \(V\)-valued map and is
    absolutely continuous as an \(L^2\)-valued map. Let
    \(v_\delta=\eta_\delta*v_{\mathrm{ev}}\), where \(\eta_\delta(t)=\delta^{-1}\eta_1(t/\delta)\) and \(\eta_1\) is a
    nonnegative smooth function with support in \([-1,1]\) and integral one.
    Uniform continuity gives \(v_\delta\to v\) uniformly into \(V\) on
    \([0,\infty)\). Moreover,
    \(v_\delta'=\eta_\delta*v_{\mathrm{ev}}'\) and
    \[
      \|v_\delta'(y)\|_{L^2}
      \leq\frac{\|\eta_1\|_\infty}{\delta}\int_{y-\delta}^{y+\delta}\|v_{\mathrm{ev}}'(t)\|_{L^2}\dd t
      \le 2\|\eta_1\|_\infty\,\mathcal M\bigl(\|v_{\mathrm{ev}}'\|_{L^2}\bigr)(y),
    \]
    where \(\mathcal Mg(y):=\sup_{r>0}(1/(2r))\int_{y-r}^{y+r}\abs{g(t)}\dd t\) is the Hardy--Littlewood maximal function on \(\mathbb R\). Since \(|y|^\alpha\in A_2(\mathbb R)\), the operator \(\mathcal M\) is bounded on \(L^2(\mathbb R,|y|^\alpha\dd y)\) \cite{Muc1972,Gra2014}. As \(\|v_{\mathrm{ev}}'\|_{L^2}\) belongs to this space, the right-hand side belongs
    to \(L^2(\mathbb R,|y|^\alpha dy)\), uniformly in \(\delta\), and vector-valued Lebesgue
    differentiation followed by dominated convergence gives
    \[
      \int_0^\infty y^\alpha\|v_\delta'-v'\|_{L^2}^2dy\longrightarrow0.
    \]
    All these maps vanish for \(y>Y+1\) when \(\delta<1\).
    Uniform \(V\) convergence on this interval also gives convergence of
    the form-energy term and of the trace in \(V\).

    For fixed \(\delta\), approximate the Bochner integral defining
    \(v_\delta\) by Riemann sums of the form \(\sum_j |I_j|\eta_\delta(y-s_j)v_{\mathrm{ev}}(s_j)\).
    They and their \(y\)-derivatives converge uniformly into \(V\), by
    uniform continuity of the integrands on a common compact set.
    Each sum is a finite sum of products \(\phi(y)h(x)\) with
    \(\phi\in C_c^\infty(\mathbb R)\) and \(h\in V\).
    For \(h_k\in C_c^\infty(\R^N\setminus\Gamma_\Sigma)\) with \(h_k\to h\) in \(V\),
    \[
      \begin{aligned}
        &\int_{\R^{N+1}_+}y^\alpha|\nabla_\Aext \bigl(\phi(y)(h_k-h)\bigr)|^2\dd x\dd y\\
        &\quad=\int_0^\infty y^\alpha\bigl(|\phi|^2q(h_k-h)+|\phi'|^2\|h_k-h\|_{L^2}^2\bigr)\dd y\longrightarrow0,
      \end{aligned}
    \]
    and the traces converge in \(V\).
    Finite sums of products \(\phi(y)h_k(x)\) belong to \(\mathcal C_\Aext\);
    their spatial supports are contained in a finite union of compact
    subsets of \(\mathbb R^N\setminus\Gamma_\Sigma\).
    Choose successively the errors of the mollification, of the Riemann sum, and of the
    approximation \(h_k\to h\) below \(1/(3n)\) in the distance
    \[
      \mathcal J[v-w]^{1/2}
      +\|v(0)-w(0)\|_V.
    \]
    The resulting diagonal sequence has both required convergences.
  \end{proof}

  \Needspace{14\baselineskip}
  \begin{lemma}[Semigroup extension and core approximation]\label{lem:extension}
    Let \(f\in C_c^\infty(\R^N\setminus\Gamma_\Sigma)\). Then
    \begin{enumerate}[label=(\roman*)]
      \item\label{it:ext-slices} For every \(y>0\), \(U_f(y)\in D(L)\) and \(\norm{U_f(y)}_{L^2}\leq\norm{f}_{L^2}\). Moreover, \(y\mapsto U_f(y)\) is continuous from \([0,\infty)\) into \(V\).
      \item\label{it:ext-energy} The map \(y\mapsto U_f(y)\) is \(C^1\) from \((0,\infty)\) into \(L^2(\R^N)\), and
        \[
          \mathcal J[U_f]=\kappa_\sigma\|L^{\sigma/2}f\|_2^2.
        \]
      \item\label{it:ext-approx} There are \(u_n\in\mathcal C_\Aext\) such that
        \[
          \mathcal J[u_n-U_f]\longrightarrow0,
        \]
        and \(u_n(\cdot,0)\to f\) in \(L^2(\R^N)\) as \(n\to\infty\).
    \end{enumerate}
  \end{lemma}

  \begin{proof}[Proof of \cref{lem:extension}]
    Since \(f\in C_c^\infty(\R^N\setminus\Gamma_\Sigma)\subset D(L)\), the proof of \cref{prop:variational-extension} applies to \(f\): for \(y>0\), \(U_f(y)=\Upsilon(yL^{1/2})f\), the map \(y\mapsto U_f(y)\) is \(C^1\) from \((0,\infty)\) into \(L^2(\R^N)\), \(U_f\in\mathcal E_f\), and the energy identity in \ref{it:ext-energy} holds. This proves \ref{it:ext-energy}. Moreover, with \(M\) and \(\phi_y\) as in that proof,
    \begin{equation}\label{eq:U-flux-bound}
      \norm{y^\alpha U_f'(y)}_{L^2}=\norm{L^\sigma\phi_y(L)f}_{L^2}\leq M\norm{L^\sigma f}_{L^2},\qquad y>0.
    \end{equation}

    \noindent\underline{\emph{Proof of \ref{it:ext-slices}.}}
    Put \(U_f(0)=f\), consistently with \(\Upsilon(0)=1\). Since \(0<\Upsilon\leq1\),
    \begin{gather*}
      \norm{U_f(y)}_{L^2}^2=\int_{[0,\infty)}\Upsilon(y\lambda^{1/2})^2\dd\mu_f\leq\norm{f}_{L^2}^2,\\
      \int_{[0,\infty)}\lambda^2\Upsilon(y\lambda^{1/2})^2\dd\mu_f\leq\int_{[0,\infty)}\lambda^2\dd\mu_f<\infty,
    \end{gather*}
    so \(U_f(y)\in D(L)\). For \(y,y'\geq0\),
    \[
      \norm{U_f(y)-U_f(y')}_V^2
      =
      \int_{[0,\infty)}(1+\lambda)\pt{\Upsilon(y\lambda^{1/2})-\Upsilon(y'\lambda^{1/2})}^2\dd\mu_f(\lambda),
    \]
    which tends to zero as \(y'\to y\) by dominated convergence, because \(\int(1+\lambda)\dd\mu_f=\norm{f}_V^2<\infty\).\smallskip

    \noindent\underline{\emph{Proof of \ref{it:ext-approx}.}} Let \(\theta\) be as in the proof of \cref{prop:Kstar-le-Kred}, and for \(R\geq1\) let \(v_R(y):=\theta(y/R)U_f(y)\). By \ref{it:ext-slices}, \(v_R\) is continuous from \([0,\infty)\) into \(V\), \(v_R(0)=f\), and \(v_R(y)=0\) for \(y\geq2R\). By \ref{it:ext-energy}, \(v_R\) is \(C^1\) on \((0,2R]\) into \(L^2(\R^N)\), and by \ref{it:ext-slices} and \eqref{eq:U-flux-bound}, \(\norm{v_R'(y)}_{L^2}\leq C\pt{\norm{f}_{L^2}+y^{-\alpha}\norm{L^\sigma f}_{L^2}}\), which is integrable near \(0\) because \(\alpha<1\). Hence \(v_R\) is absolutely continuous on \([0,2R]\) into \(L^2(\R^N)\), and \(\int_0^{2R}y^\alpha\norm{v_R'}_{L^2}^2\dd y<\infty\) because \(U_f\) has finite energy. Since \(\theta(y/R)\) does not depend on \(x\), \(q\pt{v_R(y)-U_f(y)}=\pt{1-\theta(y/R)}^2q\pt{U_f(y)}\), and
    \[
      \begin{aligned}
        &\mathcal J[v_R-U_f]\\
        &\quad\leq
        2\int_0^\infty y^\alpha\pt{1-\theta\pt{\tfrac yR}}^2\pt{q\pt{U_f(y)}+\norm{U_f'(y)}_{L^2}^2}\dd y+
        \frac{2\norm{\theta'}_\infty^2}{R^2}\int_R^{2R}y^\alpha\norm{U_f(y)}_{L^2}^2\dd y.
      \end{aligned}
    \]
    The first term tends to zero as \(R\to\infty\) by dominated convergence and \ref{it:ext-energy}. By \ref{it:ext-slices}, the second is at most \(2\norm{\theta'}_\infty^2\norm{f}_{L^2}^2(2^{\alpha+1}/(\alpha+1))R^{\alpha-1}\), which tends to zero because \(\alpha<1\). By \cref{lem:approximation}, for each \(R\) there is \(u\in\mathcal C_\Aext\) such that
    \[
      \mathcal J[u-v_R]\leq\frac1R
    \]
    and \(\norm{u(\cdot,0)-f}_{L^2}\leq1/R\). By the triangle inequality for \(\mathcal J[\cdot]^{1/2}\), these functions satisfy~\ref{it:ext-approx} as \(R\to\infty\).
  \end{proof}

  \subsection{The proof of \cref{thm:fractional}}
  First take \(f\in C_c^\infty(\R^N\setminus\Gamma_\Sigma)\).
  By \cref{lem:extension}, there are \(u_n\in\mathcal C_\Aext\) that converge to \(U_f\) in
  energy and their traces converge to \(f\) in \(L^2\). Applying
  \eqref{eq:trace-C} to differences shows that these traces also converge in
  \(L^2(|x|^{-2\sigma}\dd x)\), with the same limit. Consequently,
  \[
    K_{N,\a,\alpha}\int_{\R^N}|x|^{-2\sigma}|f|^2\dd x\leq\mathcal J[U_f]=\kappa_\sigma\|L^{\sigma/2}f\|_2^2,
  \]
  and \(\mathcal S_{N,\sigma,\a}\geq K_{N,\a,\alpha}/\kappa_\sigma\).
  Conversely, for \(u\in\mathcal C_\Aext\) with \(u(\cdot,0)=f\), \cref{prop:variational-extension} gives \(\mathcal J[U_f,u-U_f]=0\). Hence
  \[
    \begin{aligned}
      \int_{\R^{N+1}_+}y^\alpha|\nabla_\Aext u|^2\dd x\dd y
      &\;=\kappa_\sigma\|L^{\sigma/2}f\|_2^2+\mathcal J[u-U_f]\\
      &\;\geq\kappa_\sigma\|L^{\sigma/2}f\|_2^2\\
      &\;\geq\kappa_\sigma\,\mathcal S_{N,\sigma,\a}\int_{\R^N}|x|^{-2\sigma}|f|^2\dd x.
    \end{aligned}
  \]
  Taking the trace infimum proves \eqref{eq:S-equals-K}. Substitution of
  \(\alpha=1-2\sigma\) and the Gamma recurrence yields
  \[
    \frac{K_{N,\a,\alpha}}{\kappa_\sigma}
    =2^{2\sigma}\frac{\Gamma((\mu+2+2\sigma)/4)^2}
    {\Gamma((\mu+2-2\sigma)/4)^2}.
  \]
  Since \(\mathcal S_{N,\sigma,\a}\) is an infimum over \(C_c^\infty(\R^N\setminus\Gamma_\Sigma)\), this proves \eqref{eq:fractional-inequality} on this space, with the best constant.

  We now extend \eqref{eq:fractional-inequality} to the form domain \(D(L^{\sigma/2})\) of \(L^\sigma\). Let \(f\in D(L^{\sigma/2})\), and write \(\norm{g}_\sigma^2:=\norm{g}_2^2+\|L^{\sigma/2}g\|_2^2\) for the graph norm of \(L^{\sigma/2}\). For \(g\in V=D(L^{1/2})\), the inequality \(\lambda^\sigma\leq1+\lambda\) gives
  \begin{equation}\label{eq:graph-V}
    \norm{g}_\sigma^2=\int_{[0,\infty)}(1+\lambda^\sigma)\dd\mu_g(\lambda)\leq\int_{[0,\infty)}(2+\lambda)\dd\mu_g(\lambda)\leq2\norm{g}_V^2.
  \end{equation}
  We construct \(f_n\in C_c^\infty(\R^N\setminus\Gamma_\Sigma)\) with \(\norm{f_n-f}_\sigma\to0\).
  \begin{enumerate}[label=$\bullet$,itemsep=0em]
    \item For \(M>0\), let \(f_M:=\mathbf1_{[0,M]}(L)f\). Then \(f_M\in V\), because \(\int_{[0,M]}\lambda\dd\mu_f\leq M\norm{f}_2^2\), and \(\norm{f-f_M}_\sigma^2=\int_{(M,\infty)}(1+\lambda^\sigma)\dd\mu_f\to0\) as \(M\to\infty\), because \(\int_{[0,\infty)}(1+\lambda^\sigma)\dd\mu_f<\infty\).
    \item Since \(V\) is the closure of \(C_c^\infty(\R^N\setminus\Gamma_\Sigma)\) in \(\norm{\cdot}_V\), for every \(n\) we may choose \(M_n\) with \(\norm{f-f_{M_n}}_\sigma<1/n\), and then \(f_n\in C_c^\infty(\R^N\setminus\Gamma_\Sigma)\) with \(\norm{f_n-f_{M_n}}_V<1/n\). By \eqref{eq:graph-V}, \(\norm{f_n-f}_\sigma<(1+\sqrt2)/n\).
  \end{enumerate}
  By \eqref{eq:fractional-inequality} on \(C_c^\infty(\R^N\setminus\Gamma_\Sigma)\), applied to \(f_n-f_k\), the sequence \((f_n)\) is a Cauchy sequence in \(L^2(\R^N,\abs x^{-2\sigma}\dd x)\); let \(g\) be its limit. Since also \(f_n\to f\) in \(L^2(\R^N)\), a subsequence converges almost everywhere both to \(g\) and to \(f\), so \(g=f\). Letting \(n\to\infty\) in \eqref{eq:fractional-inequality} for \(f_n\) proves \eqref{eq:fractional-inequality} for \(f\). Since \(C_c^\infty(\R^N\setminus\Gamma_\Sigma)\subset D(L^{\sigma/2})\), the best constant on \(D(L^{\sigma/2})\) is therefore also \(\mathcal S_{N,\sigma,\a}\).

  It remains to show that equality in \eqref{eq:fractional-inequality} for \(f\) forces \(f=0\).
  Let \(f_n\in C_c^\infty(\R^N\setminus\Gamma_\Sigma)\) be the approximants just constructed and let
  \(U_n:=U_{f_n}\). Each \(U_n\) belongs to \(\mathcal D_\Aext\)
  by \cref{lem:extension,prop:energy-norm}. Linearity and the energy identity give
  \[
    \mathcal J[U_n-U_m]=\kappa_\sigma\|L^{\sigma/2}(f_n-f_m)\|_2^2.
  \]
  Thus \(U_n\) converges in the homogeneous completion to a function
  \(U\) with weighted trace \(f\) and
  \[
    \int_{\R^{N+1}_+}y^\alpha|\nabla_\Aext U|^2\dd x\dd y=\kappa_\sigma\|L^{\sigma/2}f\|_2^2.
  \]
  Equality in the fractional inequality would give
  \[
    \int_{\R^{N+1}_+}y^\alpha|\nabla_\Aext U|^2\dd x\dd y=K_{N,\a,\alpha}\int_{\R^N}|x|^{\alpha-1}|U(x,0)|^2\dd x.
  \]
  By \cref{prop:energy-norm}, \(U=0\), and therefore \(f=0\).
  \qed

  \section{Examples}\label{sec:examples}
  In this section we compute \(\lam\) for Aharonov--Bohm potentials in dimensions two and three and for a Dirac monopole. By \cref{thm:main,thm:fractional}, the sharp constants depend on \(\a\) and \(\Sigma\) only through \(\lam\).

  \subsection{The case \texorpdfstring{\(N=2\)}{N=2}}\label{sec:example-2D}
  In this subsection we compute \(\lam\) for Aharonov--Bohm potentials on the full circle and on a circle with one point excluded.

  We take \(\Sigma=\emptyset\) and use polar coordinates to write \(x=r(\cos\theta,\sin\theta)\) for \(r>0\) and \(\theta\in [0,2\pi)\). In this case, an Aharonov--Bohm type vector potential takes the form of a smooth vector field \(\A:\R^2\setminus\{0\}\to \R^2\) satisfying
  \[
    \A(r,\theta):=\frac{1}{r}a_\theta(\theta) (-\sin\theta,\cos\theta)
  \]
  for some smooth real-valued \(2\pi\)-periodic function \(a_\theta\). Thus the tangential vector field of the introduction is \(\a(\theta)=a_\theta(\theta)e_\theta\), where \(e_\theta=(-\sin\theta,\cos\theta)\). The corresponding magnetic Laplace operator on the circle can be written as
  \[
    \Delta_\a u
    =
    (\partial_\theta-\ii a_\theta(\theta))^2 u
    =
    \partial_\theta^2 u-2\ii a_\theta(\theta)\partial_\theta u-\ii a_\theta'(\theta)u-a_\theta(\theta)^2u.
  \]

  Define the normalized flux and the gauge phase by
  \[
    \bar\a:=\frac{1}{2\pi}\int_0^{2\pi}a_\theta(s)\dd s,
    \qquad
    \varphi(\theta):=\int_0^\theta a_\theta(s)\dd s.
  \]
  Writing \(u=e^{\ii\varphi}u_{\mathrm g}\) gives
  \[
    (\partial_\theta-\ii a_\theta)u=e^{\ii\varphi}u_{\mathrm g}',
    \qquad
    -\Delta_\a u=e^{\ii\varphi}(-u_{\mathrm g}'').
  \]
  Thus the angular problem is the ordinary second derivative with the
  twisted boundary condition
  \(u_{\mathrm g}(\theta+2\pi)=e^{-2\pi\ii\bar\a}u_{\mathrm g}(\theta)\).
  Its eigenvalues and normalized eigenfunctions are
  \[
    \lambda_m=(m-\bar\a)^2,\qquad
    u_{\mathrm g}^{(m)}(\theta)=\frac{e^{\ii(m-\bar\a)\theta}}{\sqrt{2\pi}},
    \qquad m\in\Z.
  \]
  Multiplying by \(e^{\ii\varphi}\) returns the magnetic eigenfunctions.
  If \(m_0\) is a nearest integer to \(\bar\a\), then
  \[
    \lambda_{1,\a}=\dist(\bar\a,\Z)^2,\qquad
    \psi_1(\theta)=\frac{e^{\ii(m_0-\bar\a)\theta+\ii\varphi(\theta)}}{\sqrt{2\pi}}.
  \]
  In particular, \(|\psi_1|=1/\sqrt{2\pi}\), and \eqref{hyp:lambda} holds exactly for
  nonintegral flux; then \cref{thm:main,thm:fractional} apply, the former for \(0\leq\beta_2<1\). Integer flux gives \(\lam=0\); this is the planar case excluded by \eqref{hyp:lambda}, see \cref{rem:range}.

  \medskip
  \noindent\emph{A circle with one point excluded.}
  Let
  \(\Sigma=\{(1,0)\}\), so that its complement is the interval
  \(0<\theta<2\pi\). For any real \(a_\theta\in C^\infty(0,2\pi)\),
  fix \(\theta_0\) in this interval and set
  \[
    g(\theta)=\exp\left(\ii\int_{\theta_0}^{\theta}a_\theta(s)\dd s\right).
  \]
  Multiplication by \(g\) maps \(C_c^\infty(0,2\pi)\) onto itself and satisfies \((\partial_\theta-\ii a_\theta)(gv)=gv'\).
  It therefore extends to an isometry from \(H_0^1(0,2\pi)\) onto \(H^1_\a(\Ss^1)\), even if the potential is unbounded at the endpoints. In particular,
  \[
    H^1_\a(\Ss^1)=gH_0^1(0,2\pi),\qquad \lambda_n=\frac{n^2}{4}\quad(n\geq1),\qquad \psi_1(\theta)=\frac1{\sqrt\pi}\,g(\theta)\sin(\theta/2).
  \]
  Thus \(\lam=\lambda_{1,0}^{\Sigma}=1/4\) for every real \(a_\theta\in C^\infty(0,2\pi)\). In particular, \eqref{hyp:lambda} holds, and \cref{thm:main,thm:fractional} apply, the former for \(0\leq\beta_2<1\).

  \subsection{The 3D Aharonov--Bohm case}\label{sec:example-AB3D}
  Assume \(N=3\), and let
  \[
    \A(r,\vartheta,\varphi)
    =
    \frac{\bar\a}{r\sin\vartheta} e_\varphi,
    \qquad
    \bar\a\in\R,
  \]
  where \((r,\vartheta,\varphi)\) are the usual spherical coordinates in \(\R^3\setminus\{x_1=x_2=0\}\). Then
  \[
    \A(r,\vartheta,\varphi)=\frac1r \a(\vartheta,\varphi),
    \qquad
    \a(\vartheta,\varphi)=\frac{\bar\a}{\sin\vartheta} e_\varphi.
  \]
  In Cartesian coordinates this potential is
  \[
    \A(x)
    =
    \frac{\bar\a\,(-x_2,\,x_1,\,0)}{x_1^2+x_2^2}
    =
    \pt{\A^{2D}(x_1,x_2),\,0},
  \]
  where \(\A^{2D}\) is the planar Aharonov--Bohm potential \eqref{eq:AB-potential}. Observe that this potential is singular on the whole \(x_3\)-axis, it is homogeneous of degree \(-1\), and transverse, \(\A(x)\cdot x=0\). The field \(\a\) is smooth on \(\Ss^2\) except at the two poles, so in the framework of the introduction \(\Sigma\) consists of the two poles and \(\Gamma_\Sigma\) is the \(x_3\)-axis.

  Let \(m_0\in\Z\) satisfy \(|\bar\a-m_0|=\dist(\bar\a,\Z)\), set \(b:=\dist(\bar\a,\Z)\), and choose the normalization constant \(c_b>0\) so that
  \[
    \psi_1(\vartheta,\varphi)
    :=
    c_b (\sin\vartheta)^{b} e^{\ii m_0\varphi}
  \]
  has \(L^2(\Ss^2)\)-norm one. We show that \(\lam=b(b+1)\) and that \(\psi_1\) is a first eigenfunction. In the coordinates \(z=\cos\vartheta\), the magnetic energy \(\int_{\Ss^2}|\nabla_{\a}\psi|^2\dd\omega\) is the Dirichlet form of \cite[(2.11)]{BDELL2021} with flux \(a=\bar\a\), and \cite[Lemma 2.2]{BDELL2021} gives
  \[
    \int_{\Ss^2}|\nabla_{\a}\psi|^2\dd\omega\geq b(b+1)\int_{\Ss^2}|\psi|^2\dd\omega,
    \quad\text{ for }\; \psi\in C_c^\infty(\Ss^2\setminus\Sigma).
  \]
  By closure, \(\lam\geq b(b+1)\). Conversely, \(\psi_1\) is bounded and smooth away from the poles, and
  \[
    |\nabla_{\a}\psi_1|^2=|\partial_\vartheta\psi_1|^2+\frac{|\partial_\varphi\psi_1-\ii\bar\a\psi_1|^2}{\sin^2\vartheta}
    =c_b^2b^2(\sin\vartheta)^{2b-2}(1+\cos^2\vartheta).
  \]
  If \(b=0\), the magnetic gradient of \(\psi_1\) vanishes. If \(b>0\), let \(I_n:=\int_0^\pi(\sin\vartheta)^n\dd\vartheta\), which is finite for \(n>-1\). Since \(\dd\omega=\sin\vartheta\dd\vartheta\dd\varphi\) and \(1+\cos^2\vartheta=2-\sin^2\vartheta\),
  \[
    \int_{\Ss^2}|\nabla_{\a}\psi_1|^2\dd\omega
    =2\pi c_b^2b^2\bigl(2I_{2b-1}-I_{2b+1}\bigr),
    \qquad
    \int_{\Ss^2}|\psi_1|^2\dd\omega=2\pi c_b^2I_{2b+1},
  \]
  and both are finite because \(2b-1>-1\). The union of the two poles has zero capacity by \cref{lem:points} and a finite-sum cutoff construction, so \cref{lem:capacity-bounded} places \(\psi_1\) in \(H_{\a}^1(\Ss^2)\). For \(b=0\) its Rayleigh quotient is \(0\). For \(b>0\), integration by parts gives \(I_{2b+1}=(2b/(2b+1))I_{2b-1}\), so the Rayleigh quotient of \(\psi_1\) is
  \[
    \frac{b^2\bigl(2I_{2b-1}-I_{2b+1}\bigr)}{I_{2b+1}}
    =b^2\left(\frac{2b+1}{b}-1\right)
    =b(b+1).
  \]
  This proves the claim. Since \(N=3\), hypothesis \eqref{hyp:lambda} holds for every \(\bar\a\in\R\), and \cref{thm:main,thm:fractional} apply. Moreover, \(\lam>0\) precisely when \(\bar\a\notin\Z\).

  \begin{remark}
    The singularity of \(\a\) is essential: \(\a\) is curl-free on \(\Ss^2\setminus\Sigma\), and by \cref{prop:lambda1-zero} a smooth curl-free field on \(\Ss^2\) would have first eigenvalue \(0\).
  \end{remark}

  \subsection{The 3D Dirac monopole}\label{sec:example-monopole}
  Assume \(N=3\), and let
  \[
    \A_\nu(r,\vartheta,\varphi)
    =
    \frac{\nu(1-\cos\vartheta)}{r\sin\vartheta} e_\varphi,
    \qquad
    2\nu\in\Z.
  \]
  This is the vector potential of the Dirac monopole \cite{Dirac1931} in the north chart of Wu and Yang \cite[(5)]{WuYang1976}: away from the half-line \(\{x_1=x_2=0,\ x_3\leq0\}\), \(\nabla\times\A_\nu=\nu\,x/|x|^3\), the magnetic field of a point charge at the origin. Write
  \[
    \a_\nu(\vartheta,\varphi)
    =
    \frac{\nu (1-\cos\vartheta)}{\sin\vartheta} e_\varphi
  \]
  for the corresponding tangential field on the north chart of \(\Ss^2\).
  The excluded set is the south pole, so \(\Gamma_\Sigma=\{(0,0,x_3):x_3\leq0\}\). The monopole operator
  acts on sections of a line bundle, rather than on globally smooth
  scalar functions. Its connection one-forms on the two charts are
  \[
    \a_N\cdot\dd\omega=\nu(1-\cos\vartheta)\dd\varphi,
    \qquad \a_S\cdot\dd\omega=-\nu(1+\cos\vartheta)\dd\varphi.
  \]
  Their difference is \(2\nu\dd\varphi\). With the convention
  \(\nabla_\a=\nabla-\ii \a\), section representatives satisfy
  \(\psi_N=e^{2\ii\nu\varphi}\psi_S\); the condition \(2\nu\in\Z\) makes the transition single-valued.

  Restriction to the north chart preserves the \(L^2\)-norm and
  magnetic energy. Every function in \(C_c^\infty(\Ss^2\setminus\Sigma)\) defines a smooth global section by setting its south-chart representative equal to zero near
  the missing pole. Conversely, every smooth global section is bounded and has finite connection energy; its north-chart representative belongs to \(H^1_{\a_\nu}(\Ss^2)\) by~\cref{lem:points,lem:capacity-bounded}. Smooth sections are dense in the form domain of the smooth connection on the compact sphere. The two inclusions therefore give a unitary identification of the complete closed forms, their associated operators, and all their spectral powers. The monopole spectrum \cite[(15), (16) and (19)]{WuYang1976} (see also \cite[(1.1)]{KV2014}, where \(q=2\nu\)) is
  \[
    \ell(\ell+1)-\nu^2,\qquad
    \ell=|\nu|,|\nu|+1,\ldots,
  \]
  whose minimum is \(|\nu|\). Thus every function in \(C_c^\infty(\Ss^2\setminus\Sigma)\) has
  Rayleigh quotient at least \(|\nu|\), and \(\lambda_{1,\a_\nu}\geq|\nu|\). Conversely, let
  \[
    \psi_1(\vartheta,\varphi):=\pt{\frac{2|\nu|+1}{4\pi}}^{1/2}\pt{\cos\frac\vartheta2}^{2|\nu|}.
  \]
  It is bounded and smooth away from the south pole. Since \((1-\cos\vartheta)/\sin\vartheta=\tan(\vartheta/2)\),
  \[
    |\nabla_{\a_\nu}\psi_1|^2=|\partial_\vartheta\psi_1|^2+\nu^2\tan^2\frac\vartheta2\,|\psi_1|^2
    =\frac{2|\nu|+1}{2\pi}\,\nu^2\pt{\sin\frac\vartheta2}^2\pt{\cos\frac\vartheta2}^{4|\nu|-2}.
  \]
  If \(\nu=0\), then \(\psi_1\) is constant and has magnetic energy \(0\). If \(\nu\neq0\), the substitution \(s=\cos(\vartheta/2)\), for which \(\sin\vartheta\dd\vartheta=-4s\dd s\), gives
  \[
    \begin{aligned}
      \int_{\Ss^2}|\psi_1|^2\dd\omega&=(2|\nu|+1)\int_0^1 2s^{4|\nu|+1}\dd s=1,\\
      \int_{\Ss^2}|\nabla_{\a_\nu}\psi_1|^2\dd\omega&=4(2|\nu|+1)\nu^2\int_0^1(1-s^2)s^{4|\nu|-1}\dd s=|\nu|.
    \end{aligned}
  \]
  By \cref{lem:points,lem:capacity-bounded}, \(\psi_1\in H_{\a_\nu}^1(\Ss^2)\). Therefore \(\lambda_{1,\a_\nu}=|\nu|\), and \(\psi_1\) is a first eigenfunction.
  Since \(N=3\), hypothesis~\eqref{hyp:lambda} holds, and \cref{thm:main,thm:fractional} apply. Moreover, \(\lam>0\) precisely when \(\nu\neq0\).

  \section*{Statements and Declarations}
  \addcontentsline{toc}{section}{Statements and Declarations}

  \paragraph{Author contributions.}
  All authors contributed equally to this work.\vspace{-.7em}

  \paragraph{Data availability.}
  No datasets were generated or analyzed in this study. No external research dataset is required for these proofs.\vspace{-.7em}

  \paragraph{Conflict of interest.}
  The authors have no competing interests to declare that are relevant to the content of this article.\vspace{-.7em}

  \paragraph{Use of artificial intelligence.}
  The authors used ChatGPT (OpenAI; Astra 6.0) and Claude (Anthropic; Opus 5.5) to assist with routine computations and to improve the clarity and presentation of the manuscript. All mathematical results, computations and references were independently verified by the authors, who take full responsibility for the content of the manuscript.\vspace{-.7em}

  \paragraph{Acknowledgements.}
  AZ gratefully acknowledges financial support from the Instituto de Matem\'aticas, Universidad de Talca, through its visitors program for his research visits during the summers of 2025 and 2026. HC gratefully acknowledges financial support from the Instituto de Ciencias de la Ingenier\'ia, Universidad de O'Higgins, for his two research visits to ICI at UOH in 2025 and 2026.

  \appendix

  \section{Hypergeometric identities}\label{sec:hypergeometric}
  We collect the hypergeometric identities used in \cref{sec:general-ode}.

  We use the ordinary hypergeometric function
  \begin{equation}\label{def:hypergeometric}
    {}_2F_1(a,b;c;z)=\sum_{n=0}^\infty
    \frac{(a)_n(b)_n}{(c)_n n!}z^n,\qquad |z|<1,
  \end{equation}
  continued along the real intervals where it is used. It solves
  \begin{equation}\label{eq:gauss-ode}
    z(1-z)v''+[c-(a+b+1)z]v'-ab\,v=0.
  \end{equation}
  On \((-\infty,0)\), Kummer's solution at infinity
  \begin{equation}\label{eq:kummer-infinity}
    v(z)=(-z)^{-a}{}_2F_1\left(a,a-c+1;a-b+1;\frac1z\right)
  \end{equation}
  also solves \eqref{eq:gauss-ode}, see \cite[Vol.~I, \S2.9, (9)]{ErMagOb1953}.
  Pfaff's transformation and Gauss's summation formula are
  \begin{equation}\label{eq:pfaff}
    {}_2F_1(a,b;c;z)=(1-z)^{-a}
    {}_2F_1\left(a,c-b;c;\frac{z}{z-1}\right),
  \end{equation}
  \begin{equation}\label{eq:gauss-sum}
    {}_2F_1(a,b;c;1)=
    \frac{\Gamma(c)\Gamma(c-a-b)}{\Gamma(c-a)\Gamma(c-b)},
    \qquad \Re(c-a-b)>0.
  \end{equation}
  For \(c-a-b\notin\Z\), the connection at one reads
  \begin{equation}\label{eq:connection-one}
    \begin{aligned}
      &{}_2F_1(a,b;c;z)
      =\frac{\Gamma(c)\Gamma(c-a-b)}{\Gamma(c-a)\Gamma(c-b)}
      {}_2F_1(a,b;a+b-c+1;1-z)\\
      &\hspace{6em}+(1-z)^{c-a-b}\frac{\Gamma(c)\Gamma(a+b-c)}{\Gamma(a)\Gamma(b)}{}_2F_1(c-a,c-b;c-a-b+1;1-z).
    \end{aligned}
  \end{equation}
  These are \cite[Vol.~I, \S2.9, (3); \S2.1.3, (14); and \S2.9, (33)]{ErMagOb1953}; in the last one, the functions \(u_1\), \(u_2\) and \(u_6\) of that reference are given by \S2.9, (1), (5) and (21) there.
  For our parameters \(c-a-b=1-C\in(0,1)\), all displayed Gamma factors
  are finite and the local series in \(1-z\) justify differentiated expansions.

  \section{Sets of zero capacity}\label{sec:zero-capacity}
  In this appendix we show that points of \(\Ss^2\) have zero capacity, and that bounded functions with finite magnetic energy on the complement of a closed set of zero capacity belong to \(H_\a^1(\Sph)\).

  Following \cite[Section 13.1]{Mazya2011}, where \(\operatorname{Cap}(e,W^l_p)\) is defined through functions equal to one near \(e\), a closed set \(\Sigma\subset\Sph\) has \emph{zero capacity} if there are Lipschitz functions \(\varphi_k:\Sph\to\R\) such that each \(\varphi_k\) satisfies \(\varphi_k\geq1\) on a neighborhood of \(\Sigma\), and
  \begin{equation}\label{eq:zero-capacity}
    \int_{\Sph}\pt{\abs{\nabla_{\Sph}\varphi_k}^2+\abs{\varphi_k}^2}\dd\omega
    \xrightarrow[k\to\infty]{}0.
  \end{equation}
  A closed set of zero capacity has zero measure, because \(\abs{\Sigma}\leq\int_{\Sph}\abs{\varphi_k}^2\dd\omega\). This notion is used in \cref{prop:lambda1-zero} and in \cref{sec:example-AB3D,sec:example-monopole}.

  \begin{lemma}[Zero capacity of points on the sphere]\label{lem:points}
    A point \(p\in\Ss^2\) has zero capacity.
  \end{lemma}

  \begin{proof}[Proof of Lemma~\ref{lem:points}]
    For \(\ve\in(0,1)\), let
    \[
      \varphi_\ve(\omega):=
      \begin{cases}
        1, & \dist(\omega,p)\leq\ve^2,\\[2pt]
        \dfrac{\log(\ve/\dist(\omega,p))}{\log(1/\ve)}, & \ve^2<\dist(\omega,p)<\ve,\\[6pt]
        0, & \dist(\omega,p)\geq\ve,
      \end{cases}
    \]
    where \(\dist\) is the geodesic distance on \(\Ss^2\). Then \(\varphi_\ve\) is Lipschitz, takes values in \([0,1]\), and equals \(1\) on a neighborhood of \(p\). In geodesic polar coordinates \((r,\vartheta)\) centered at \(p\) we have \(\dd\omega=\sin r\dd r\dd\vartheta\) and \(\abs{\nabla\varphi_\ve}=(r\log(1/\ve))^{-1}\) for \(\ve^2<r<\ve\), and \(\nabla\varphi_\ve=0\) elsewhere. Since \(\sin r\leq r\),
    \[
      \int_{\Ss^2}\abs{\nabla_{\Ss^2}\varphi_\ve}^2\dd\omega
      \leq
      \frac{2\pi}{\log(1/\ve)^2}\int_{\ve^2}^{\ve}\frac{\dd r}{r}
      =
      \frac{2\pi}{\log(1/\ve)},
      \qquad
      \int_{\Ss^2}\varphi_\ve^2\dd\omega
      \leq
      2\pi\int_0^{\ve}r\dd r
      =
      \pi\ve^2,
    \]
    and both tend to zero as \(\ve\to0^+\).
  \end{proof}

  A finite union of zero-capacity closed sets also has zero capacity:
  sum the corresponding cutoff functions and truncate at one, using
  \(|\sum_{j=1}^m\nabla\varphi_j|^2
  \leq m\sum_{j=1}^m|\nabla\varphi_j|^2\).
  In particular, the union of the two poles of \(\Ss^2\) has zero capacity.

  \begin{lemma}[Removability of zero-capacity sets]\label{lem:capacity-bounded}
    Let \(\Sigma\subset\Sph\) be a closed set of zero capacity, let \(\a\) be smooth on \(\Sph\setminus\Sigma\), and let \(\psi:\Sph\setminus\Sigma\to\C\) be bounded and locally Lipschitz, with \(\int_{\Sph\setminus\Sigma}\abs{\nabla_\a\psi}^2\dd\omega<\infty\). Then \(\psi\in H_\a^1(\Sph)\). In particular, \(\int_{\Sph}\abs{\nabla_\a\psi}^2\dd\omega\geq\lam\int_{\Sph}\abs{\psi}^2\dd\omega\).
  \end{lemma}
  \begin{proof}[Proof of Lemma~\ref{lem:capacity-bounded}]
    Let \(\varphi_k\) be the cutoffs in \eqref{eq:zero-capacity} and
    put \(\chi_k=\min\{1,\max\{\varphi_k,0\}\}\). Then
    \(0\leq\chi_k\leq1\), \(\chi_k=1\) near \(\Sigma\), and
    \(\|\chi_k\|_{H^1}\to0\). The product rule is
    \(\nabla_\a(\chi_k\psi)=\chi_k\nabla_\a\psi+\psi\nabla_{\Sph}\chi_k\).
    Passing to a subsequence, we may assume \(\chi_k\to0\) almost
    everywhere, since \(\|\chi_k\|_{L^2}\to0\). Then
    \[
      \|\chi_k\nabla_\a\psi\|_{L^2}\to0,
      \qquad
      \|\psi\nabla\chi_k\|_{L^2}
      \le\|\psi\|_{L^\infty}\|\nabla\chi_k\|_{L^2}\to0,
      \qquad \|\chi_k\psi\|_{L^2}\to0.
    \]
    The first assertion follows by dominated convergence. Hence \((1-\chi_k)\psi\to\psi\) in the magnetic form norm. For each fixed \(k\) the cutoff product has compact support away from \(\Sigma\), where \(\a\) is bounded and smooth. Local smoothing in charts approximates it by \(C_c^\infty(\mathbb S^{N-1}\setminus\Sigma)\) in the same norm. A diagonal choice proves the claim.
  \end{proof}

  \smallskip
  \begingroup\small
  \bibliographystyle{abbrv}
  \bibliography{citations-Magnetic_Hardy-A-final.bib}
  \endgroup

  \medskip
  \medskip

  \begingroup
  \small
  \addressblock{Instituto de Matem\'aticas, Universidad de Talca (UTalca). Camino Lircay S/N, Campus Norte, Talca, Chile.}
  \EmailB{hcastro@utalca.cl}

  \bigskip

  \addressblock{Instituto de Ciencias de la Ingenier\'ia, Universidad de O'Higgins (UOH). Avenida del Libertador Bernardo O'Higgins 611, Rancagua, Chile.}
  \EmailB{andres.zuniga@uoh.cl}
  \endgroup

  \end{document}